\documentclass[11pt]{article}
\usepackage{mathrsfs}
\usepackage[toc,page]{appendix}
\usepackage{amsmath}
\usepackage{calrsfs}
\usepackage{mathrsfs}
\usepackage[mathscr]{euscript}
\usepackage{amsfonts}
\usepackage{amssymb, amsmath, cite}
\usepackage{color}

\newtheorem{theorem}{Theorem}
\newtheorem{lemma}{Lemma}
\newtheorem{proposition}{Proposition}
\newtheorem{remark}{Remark}

\newtheorem{corollary}{Corollary}

\newcommand\be{\begin{equation}}
\newcommand\ee{\end{equation}}
\newcommand\ber{\begin{eqnarray}}\newcommand\bea{\begin{eqnarray}}
\newcommand\eer{\end{eqnarray}}\newcommand\eea{\end{eqnarray}}
\newcommand\berr{\begin{eqnarray*}}
\newcommand\eerr{\end{eqnarray*}}

\newcommand{\lm}{\lambda}

\newcommand\re{\mathrm{e}}

\newcommand{\ud}{\mathrm{d}}
\newcommand{\nm}{\nonumber}\newcommand{\nn}{\nonumber}

\newcommand{\ito}{\int_{\Omega}}

\newcommand{\vep}{\varepsilon}

 \title{Multiple Solutions for the Non-Abelian Chern--Simons--Higgs Vortex Equations\footnote{This work is supported by PRIN12:  "Variational and Perturbative Aspects in Nonlinear Differential Problems", PRIN15: "Variational methods, with applications to problems in mathematical physics and geometry", FIRB project: "Analysis and Beyond", MIUR Excellence Project: Department of Mathematics, University of Rome Tor Vergata, CUP E83C18000100006", National Natural Science Foundation of China under Grants 11671120 and 11471100,  and  HASTIT(18HASTIT028).}}

\author{Xiaosen Han$^{a}$ \quad   \quad Gabriella Tarantello$^{b}$
   \\
\small {$^a$ Institute of Contemporary Mathematics, School of Mathematics and Statistics,  Henan University}\\
\small{Kaifeng  475004, PR China}\\
\small {$^b$ Dipartimento di Matematica, Universit\`{a} degli Studi di Roma ``Tor Vergata'',}\\
\small{Via della Ricerca Scientifica, 00133 Rome, Italy}}
\date{}

\begin{document}
\maketitle
\begin{abstract}
In this paper we study  the existence of multiple solutions for the non-Abelian Chern--Simons--Higgs $(N\times N)$-system:
\[
\Delta u_i=\lambda\left(\sum_{j=1}^N\sum_{k=1}^N K_{kj}K_{ji}\re^{u_j}\re^{u_k}-\sum_{j=1}^N K_{ji}\re^{u_j}\right)+4\pi\sum_{j=1}^{n_i}\delta_{p_{ij}},\quad i=1,\dots, N;
\]
over a doubly periodic domain $\Omega$, with coupling matrix $K$ given by the Cartan matrix of $SU(N+1),$ (see \eqref{k1} below). Here,  $\lambda>0$ is the coupling parameter, $\delta_p$ is the Dirac measure with pole at $p$ and  $n_i\in \mathbb{N},$ for  $i=1, \dots, N.$
When $N=1, 2$ many results are now available for the periodic solvability of such system and provide the existence of different classes of solutions known as: topological, non-topological, mixed and blow-up type.
On the contrary for $N\ge 3,$ only recently in  \cite{haya1} the authors managed to obtain the existence of one doubly periodic solution via a minimisation procedure, in the spirit of \cite{nota} . Our main contribution in this paper is to
show (as in \cite{nota}) that actually  the given  system admits a second doubly periodic solutions of ``Mountain-pass'' type, provided that $3\le N\le 5$. Note that the existence of multiple solutions is relevant from the physical point of view.  Indeed, it implies the co-existence  of different non-Abelian Chern--Simons condensates sharing the same  set (assigned component-wise) of vortex points, energy and fluxes.
The main difficulty to overcome is to attain  a ``compactness'' property  encompassed by the so called Palais--Smale condition  for the corresponding ``action'' functional,  whose validity remains still open for  $N\ge 6$.
\medskip

\end{abstract}

\section{Introduction:}

In recent years, the Chern--Simons forms proposed by Chern and Simons \cite{CS1,CS2}  concerning  secondary characteristic classes have played a very important role both in theoretical and applied sciences.
In this respect we mention,  knot invariants\cite{der}, Jones polynomial \cite{wt1}, quantum field theory \cite{Aty1}, string theory\cite{Mar1,wt2}, high-temperature superconductivity \cite{Kh,Ma,zT,W2}, optics \cite{Be}, and condensed matter physics \cite{In,Ray,qizh1}.

In superconductivity, Hong--Kim--Pac \cite{zHKP} and Jackiw--Weinberg \cite{JW} introduced the Chern--Simons terms into the Abelian Higgs
model  to describe particles carrying  both magnetic and electric charges. In addition, in \cite{zHKP} and \cite{JW} the authors showed that, by neglecting the Maxwell term, one could attain a self-dual BPS-regime (Bogomol'nyi\cite{zBo} and Prasad--Sommerfield\cite{zPS})  with a 6th-order potential.
Since then, many other physical Chern--Simons models have been introduced with analogous features  \cite{bhsm,Burct,Hozh,jaclw1,Kim1,Gho1}.
Starting with the work in \cite{wangr, spya4,spya3, caya1, nota1, taran96},   a rather complete description of (electro-magntic) abelian Chern--Simons vortices is now available in literature, see  \cite {tabook,Yang1} for a detailed account.

However, more recently there has been a growing interest towards non-Abelian vortices concerning particle interactions other than electro-magnetic ones (e.g. weak, strong, electro-weak  etc).
 Indeed within the general framework of Supersymmetry, it has been noted that non-Abelian vortices assume a relevant  role  towards the delicate issue of ``confinement''.  With this point of view,
 and after the ``pure'' non-Abelian Chern--Simons--Higgs model of Dunne  \cite{zD,zD2,zD3}, several other  models have been discussed in \cite{Nassfs,KK,Kaol,Oh,Leekm,Culos,Kha1},
 which have introduced also genuinely new non-Abelian ansatz  in order to attain self-duality. In this way, one can reduce the equations of motion governing non-Abelian Chern--Simons--Higgs vortices  in the (self-dual) BPS-regime into the following nonlinear elliptic system of PDE's:
\ber
\Delta u_i+\lambda\left(\sum_{j=1}^n K_{ji}\re^{u_j}-\sum_{j=1}^n\sum_{k=1}^nK_{kj}K_{ji}\re^{u_j}\re^{u_k}\right)=4\pi\sum_{j=1}^{n_i}\delta_{p_{ij}},\quad i=1,\cdots,n,\label{i1}
\eer
with a suitable coupling matrix  $K=(K_{ij})$ determined by the physical model under consideration. In  \eqref{i1}, we have  $\lambda>0$ a coupling parameter and  $n_i\in \mathbb{N}$ is the number of  assigned (vortex) points $p_{i1},\dots, p_{in_i} $  (counted with multiplicity)  for the  $i$-th component,  $ i=1, \dots, n.$

For the ``pure''  Chern--Simons--Higgs  model in  \cite{zD}, the matrix $K=(K_{ij})$ coincides  with  the Cartan matrix corresponding to the (non-Abelian) gauge group $G$ describing the internal symmetries of the model. Typically $G$ admits a finite-dimensional semi-simple Lie algebra $L$ and $n=\text{ rank } L. $\\
For example, to the gauge group $G=SU(N+1)$  of rank $N$ corresponds the following $N\times N$  Cartan matrix $K:$
\be
K\equiv\begin{pmatrix}
 2 & -1 & 0&\dots&\dots & 0\\
-1 &  2 & -1&0&\dots & 0\\
0 &-1&  2 &-1&\dots&0\\
\vdots& &\ddots&\ddots&\ddots&\vdots\\
0&\quad &\ddots&-1&2&-1\\
 0 & \dots &   &0&-1&2
\end{pmatrix}. \label{k1}
\ee
The first rigorous existence result about the system  \eqref{i1} in $\mathbb{R}^2$  is due to Yang  \cite{Yang2},  who uses a direct minimisation approach to establishes a planar (topological) solution
for  a general  class of coupling matrices $K,$ which  include all possible choices of Cartan matrices.  The existence of non-topological planar solutions was pursued by a perturbation approach (in the spirit of \cite{chim}) for the  Lie-Algebras of rank 2 given respectively by  $A_2$, $B_2$ and $G_2$  in \cite{alw1,alw2}, see also \cite{chkl1}. While, the existence of mixed-type planar $SU(3)$-vortices can be found in \cite{chkll1,chkl2}. See also  \cite{Huli1,Hull} for results in the skew-symmetric case.

The periodic case was first dealt in \cite{nota}, where the authors proved the existence of  multiple doubly periodic $SU(3)$-vortices, solution of \eqref{i1}--\eqref{k1} with $N=2$. This result was extended in \cite {hata}, where  \eqref{i1} is considered with a general $2\times2$ nonsingular coupling matrix $K,$ including all Cartan's type.
 See also  \cite{linya1} for the construction of bubbling solutions.

However, when the system \eqref{i1} involves three or more  components over a doubly periodic domain, then the results available are less satisfactory.  In fact, only recently Han--Yang \cite{haya1} were able to extend the constraint-minimisation approach of  \cite {nota} and established the existence of a doubly periodic solution for the system \eqref{i1}--\eqref{k1} with $N\ge3.$ A possible extension of \cite{haya1} to the system \eqref{i1} with a more general  $n\times n$ nonsingular  coupling matrix  $K$ of Cartan-type, was claimed in  Han--Lin--Yang \cite{haly}.

The aim of this paper is to show that actually  \eqref{i1}--\eqref{k1} with $N\ge3$ admits a second doubly periodic solution (other than the topological one in \cite{haya1}), which we obtain via a min-max procedure of  ``mountain-pass'' type \cite {amra}.
As already mentioned, the multiple solvability of the system \eqref{i1} is relevant from the physical point of view. Indeed, it indicates  that  (asymptotically)  each ``vacua"  states of the system may   support  a vortex configuration with the same set of vortex  points (assigned component-wise) at the
same  (qunatized) energy level. The main difficulty to apply a variational approach to the ``action'' functional corresponding to \eqref{i1} (see  \eqref{a20} below) is to show that it satisfies a ``compactness'' property, expressed by the so called Palais--Smale (PS)-condition. Such condition becomes rather involved when we deal with  three or more components, which might allow enough room for a compactness loss. We manage to resolve such a compactness issue for  \eqref{i1}--\eqref{k1}, when $N=3,4,5$, and prove the following:




\begin{theorem}\label{th1}
Consider the non-Abelian Chern--Simons--Higgs system \eqref{i1} over a doubly periodic domain $\Omega$ and with the matrix $K$  in  \eqref{k1} (i.e the Cartan matrix of $SU(N+1)$).
For $N=3,4,5$  and any given set of points   $p_{j1}, \dots, p_{jn_j}$ $ (j=1, \dots, N)$ on $\Omega$ repeated  with multiplicity,
there exists a large  constant $\lambda_1>0$ such that when $\lambda>\lambda_1$ the system \eqref{i1} admits at least two distinct  solutions.
\end{theorem}

\begin{remark}
  The constant $\lambda_1$ in our statement satisfies the following  lower bound:
 \ber
  \lambda_1>\lambda_0\equiv\frac{16\pi}{|\Omega|}\frac{\sum\limits_{i=1}^N\sum\limits_{j=1}^N(K^{-1})_{ij}n_j}{\sum\limits_{i=1}^N\sum\limits_{j=1}^N(K^{-1})_{ij}}.
 \eer
In fact the condition: $\lambda>\lambda_0$ is necessary for the existence of a doubly periodic  solution of \eqref{i1}-\eqref{k1}, as shown in \cite{haya1}.

\end{remark}

The rest of our paper is organised as follows.  In Section \ref{sct2} we present the variational formulation of the problem and furnish a new approach (different from \cite{haya1}) to solve
the associated constraint equations for the system \eqref{i1}--\eqref{k1}.  In Section \ref{sct3} we prove our main theorem by showing first that the solution obtained in
 \cite{haya1} corresponds  to a local minimum  for the ``action'' functional $I$ in \eqref{a20} below, which we show then to admits a mountain-pass structure \cite{amra}. Section \ref{sct4}
is devoted to the proof of the Palais--Smale-condition. The last section is a linear algebra  Appendix which contains useful facts needed in Section \ref{sct2}.

\section{Variational formulation and resolution of the constraints }\label{sct2}
\setcounter{equation}{0}\setcounter{lemma}{0}\setcounter{proposition}{0}\setcounter{remark}{0}\setcounter{theorem}{0}

In this section, we  carry out a variational  formulation for  \eqref{i1} and  solve the associated constrained problem  when $K$ is the Cartan matrix \eqref{k1} of $SU(N+1)$,  with $N\ge3$. It is well known that $K$ in \eqref{k1} is non-degenerate and positive definite.

Moreover by setting:
 \ber
  A \equiv K^{-1}=(a_{ij})_{N\times N},\label{a6a}
 \eer
we easily check that,
 \be
 a_{jk}=a_{kj}=\frac{1}{N+1}\Big[\min\{j, k\}(N+1-\max\{j, k\})\Big], \quad j,k=1, \dots N.\label{re0}
 \ee
Let
 \be
 r_j\equiv\sum\limits_{k=1}^Na_{jk}=\frac12j(N+1-j), \quad j=1, \dots, N \label{re}
 \ee
 and note that,
 \[
  \sum\limits_{j=1}^Nr_j=\frac{N(N+1)(N+2)}{12}. \label{2.3bis}
 \]
 Furthermore, consistently with \eqref{re},   it is convenient to set
 \be r_j=0,\, \text{for}\, \,j\le 0 \,\,\text{or}\,\, j\ge N+1.\label{re2}
 \ee

 Define,
 \be
 R\equiv{\rm diag}\big\{r_1, \dots,  r_N\big\}. \label{re3}
 \ee
 and let
 \be
 M\equiv RKR=\begin{pmatrix}
 2r_1^2 & -r_1r_2 & 0&\dots&\dots & 0\\
-r_1r_2&  2r_2^2 & -r_2r_3&0&\dots & 0\\
0 &-r_2r_3&  2r_3^2 & -r_3r_4&\dots&0\\
\vdots& &\ddots&\ddots&\ddots&\vdots\\
0& &\ddots&-r_{N-2}r_{N-1}&2r_{N-1}^2&-r_{N-1}r_N\\
 0 & \dots &   &0&-r_{N-1}r_N&2r_N^2
\end{pmatrix}.\label{me}
\ee

In what follows we replace the given unknown $u_i$ by its   translation $u_i\to u_i+\ln r_i$, which (by an abuse of notation) we still denote by $u_i$, namely:
\be
u_i\to u_i+\ln r_i, \quad   i=1, \dots, N, \label{a7}
\ee
with  $r_i$   given by \eqref{re}.

 Furthermore  we use the following  notations:
\ber
  &&\mathbf{u}\equiv(u_1, \dots, u_N)^\tau, \,  \mathrm{U}\equiv{\rm diag}\big\{\re^{u_1},  \dots, \re^{u_N}\big\},\, \mathbf{U}\equiv(\re^{u_1},  \dots, \re^{u_N})^\tau, \label{a11}\\
&& \mathbf{1}\equiv(1, \dots, 1)^\tau,\quad \mathbf{s}\equiv\left(\sum\limits_{s=1}^{n_1}\delta_{p_{1s}},  \dots,   \sum\limits_{s=1}^{n_N}\delta_{p_{Ns}}\right)^\tau,\label{a12}
\eer
which help us to  write  \eqref{i1} as follows:
\ber
 \Delta \mathbf{u}=\lambda K\mathrm{U}M(\mathbf{U}-\mathbf{1}) +4\pi\mathbf{s}, \label{a9}
\eer
once we take into account that,
 \be M\mathbf{1 }= R\mathbf{1}. \label{2.9bis}\ee

To find a doubly periodic solution of  \eqref{i1},  we  define the following  background functions\cite{aubi},
\ber
 \Delta u_i^0=4\pi\sum\limits_{s=1}^{n_i}\delta_{p_{is}}-\frac{4\pi n_i}{|\Omega|},\quad \ito u_i^0\ud x=0, \label{a.2.11}
\eer
and observe that $\re^{u_i^0}\in L^\infty(\Omega), \,  \forall i=1, \dots, N$. We set $u_i=u_i^0+v_i, \, i=1,  \dots,  N$, and we will use the following  $N$-vector  notation:
\be
 \mathbf{v}\equiv(v_1, \dots, v_N)^\tau, \quad \mathbf{n}\equiv(n_1, \dots, n_N)^\tau, \quad \mathbf{0}\equiv(0, \dots, 0)^\tau. \label{a13}
 \ee
In this way,  the  system   \eqref{a9} can be rewritten component-wise as follows:
 \be
\Delta v_i=\lm\left(\sum_{j=1}^N\sum_{k=1}^NK_{jk}K_{ij}r_jr_k\re^{u_j^0+v_j}\re^{u_k^0+v_k}-\sum_{j=1}^NK_{ij}r_j\re^{u_j^0+v_j}\right)+\frac{4\pi n_i}{|\Omega|},\quad i=1,\dots,N.\label{a14}
\ee

 To formulate the system \eqref{a14}  in a variational form,
 as in\cite{haly}, we  rewrite   \eqref{a14}  equivalently as follows:
\ber
 \Delta A \mathbf{v}&=&\lambda\mathrm{U}M(\mathbf{U}-\mathbf{1})+ \frac{\mathbf{b}}{|\Omega|}, \label{a17}
\eer
where, the matrices  $A$ and $M$ are defined  in  \eqref{a6a} and \eqref{me}  respectively,  and we have set:
 \ber
  \mathbf{b}\equiv(b_1, \dots, b_N)^\tau \quad \text{with}\quad b_j=4\pi \sum\limits_{k=1}^Na_{jk}n_k>0, \,j=1, \dots, N. \label{a18}
 \eer

Since  the matrices $A$ and $M$ defined in \eqref{a6a} are symmetric, we   obtain a variational formulation for   the system \eqref{a17},
by considering the following  (action) functional:
 \ber
  I(\mathbf{v})&=&\frac12\sum\limits_{i=1}^2\ito\partial_i\mathbf{v}^\tau A\partial_i\mathbf{v}\ud x+\frac\lambda2\ito(\mathbf{U}-\mathbf{1})^\tau M(\mathbf{U}-\mathbf{1}) \ud x
  +\frac{1}{|\Omega|}\ito \mathbf{b}^\tau\mathbf{v}\ud x. \label{a20}
 \eer

Indeed, the functional \eqref{a20} is well-defined and of class $C^1$ on the Hilbert  (product)  space  $(W^{1,2}(\Omega))^N$ considered with  the usual norm:
 \berr
 \|\mathbf{w}\|^2=\|\mathbf{w}\|_2^2+\|\nabla \mathbf{w}\|_2^2=\sum\limits_{i=1}^N\ito(w_i^2+|\nabla w_i|^2)\ud x,
 \eerr
 for any $\mathbf{w}=(w_1,\dots, w_N)^\tau, \, w_i\in W^{1,2}(\Omega), \, i=1, \dots, N.$

 It is easy to check that every critical point of $I$ in $(W^{1,2}(\Omega))^N$ defines a (weak) solution for \eqref{a17}. Although $I$ is not bounded from below,
 we show that it admits a local minimum.

 To this purpose, it is useful to consider a constrained minimization problem, firstly introduced in \cite{caya1} for the abelian Chern--Simons--Higgs equation, and
 subsequently refined in \cite{nota, haya1}  for the non-Abelian Chern--Simons--Higgs system \eqref{i1}.  The main difficulty to pursue such a constraint approach
 is to show that the given ``natural'' constraints are actually uniquely solvable with respect to the mean value of each component.

 To be more precise, we use the decomposition: $W^{1,2}(\Omega)=\mathbb{R}\oplus  \dot{W}^{1,2}(\Omega), $
  where,
\[\dot{W}^{1,2}(\Omega)\equiv\left\{w\in W^{1,2}(\Omega)\Bigg| \ito w\ud x=0\right\}\]
  is a closed subspace of $W^{1,2}(\Omega)$.
 Therefore, for any $v_i\in W^{1, 2}(\Omega)$, we set $v_i=c_i+w_i$,  with $w_i\in  \dot{W}^{1,2}(\Omega),$ and $ c_i=\frac{1}{|\Omega|}\ito v_i\ud x$,  $ i=1, \dots, N. $
 Consequently, the integration of \eqref{a14} over $\Omega$ gives the following natural constraints:
 \ber
 &&2r_j\re^{2c_j}\ito\re^{2(u_j^0+w_j)}\ud x-\re^{c_j}\ito\re^{u_j^0+w_j}\left(1+r_{j-1}\re^{c_{j-1}}\re^{u_{j-1}^0+w_{j-1}}+r_{j+1}\re^{c_{j+1}}\re^{u_{j+1}^0+w_{j+1}}\right)\ud x\nn\\
 &&+\frac{b_j}{\lambda r_j}=0, \,  j=1, \dots, N. \label{a26}
 \eer

Clearly, for any $\mathbf{w}=(w_1, \dots,w_N)^\tau$ with $w_i\in \dot{W}^{1, 2}(\Omega) \,  i=1, \dots, N$,  the equations \eqref{a26}  are solvable  with respect to $\re^{c_j}$  only if,
 \ber
&& \left(\ito\re^{u_j^0+w_j}\left(1+r_{j-1}\re^{c_{j-1}}\re^{u_{j-1}^0+w_{j-1}}+r_{j+1}\re^{c_{j+1}}\re^{u_{j+1}^0+w_{j+1}}\right)\ud x\right)^2\ge \frac{8b_j}{\lambda}\ito\re^{2(u_j^0+w_j)}\ud x,\nn\\
 && \, j=1, \dots, N. \label{a28}
 \eer
 On the other hand,  \eqref{a28} can be ensured by  requiring that,
 \ber
\left(\ito\re^{u_j^0+w_j}\ud x\right)^2&\ge& \frac{8b_j\ito\re^{2(u_j^0+w_j)}\ud x}{\lambda}, \quad j=1, \dots, N.\label{a30}
 \eer

Thus,  we  define the admissible set:
\ber
\mathcal{A}&\equiv&\left\{(w_1, \dots, w_N)^\tau \, \Big|w_j\in \dot{W}^{1, 2}(\Omega)\,\text{satisfies} \,\eqref{a30},  \,\forall j=1, \dots, N\right\}.\label{a.0}
\eer

Therefore,  for any $\mathbf{w}\in \mathcal{A}$,  to get a solution of  \eqref{a26}, it is equivalent to  show that,
\ber
\re^{c_j}=\frac{1}{4r_j\ito\re^{2(u_j^0+w_j)}\ud x}\left\{Q_j-(-1)^{\vep_j}\sqrt{Q_j^2-\frac{8b_j\ito\re^{2(u_j^0+w_j)}\ud x}{\lambda}}\right\},\, j=1, \dots, N,\label{a219}
\eer
admits a (unique) solution, with fixed $\vep_j\in\{0,1\}$ and,
\ber
Q_j\equiv\ito\re^{u_j^0+w_j}\left(1+r_{j-1}\re^{c_{j-1}}\re^{u_{j-1}^0+w_{j-1}}+r_{j+1}\re^{c_{j+1}}\re^{u_{j+1}^0+w_{j+1}}\right)\ud x.\label{a220}
\eer

In \cite{haly} the above equations \eqref{a219} are shown to be uniquely solvable when one takes $\vep_j=1, \forall j=1, \dots, N$.
In  what follows, we shall handle such a uniqueness solvability issue of \eqref{a219}, for any choice of  $\vep_j\in\{0,1\}$.

To this purpose we set $t_j=\re^{c_j}>0, \,  j=1,\dots N$ and we  show that,  for any assigned $\vep_j\in\{0,1\}$, the $N$-system of equations:
 \ber
 t_j-\frac{1}{4r_j\ito\re^{2(u_j^0+w_j)}\ud x}\left[\hat{Q}_j-(-1)^{\vep_j}\sqrt{\hat{Q}_j^2-\frac{8b_j}{\lambda}\ito\re^{2(u_j^0+w_j)}\ud x}\right]=0,\, j=1, \dots, N,\label{a.1}
  \eer
admits a unique non-degenerate solution,  smoothly depending on $(w_1, \dots, w_N)$, where
\ber
\hat{Q}_j\equiv\ito\re^{u_j^0+w_j}\left(1+r_{j-1}t_{j-1}\re^{u_{j-1}^0+w_{j-1}}+r_{j+1}t_{j+1}\re^{u_{j+1}^0+w_{j+1}}\right)\ud x.\label{a.1'}
\eer

For fixed $\vep_j\in\{0, 1\}, j=1, \dots, N$ we set:
\ber
\vep=(\vep_1, \dots, \vep_N), \label{a.2}
\eer
and for $s\in[0,1]$ we consider the following one-parameter family of functions
\ber
&&\phi_{j,s}(t_{j-1},t_j, t_{j+1}, \vep_j)\nn\\
&&\equiv t_j-\frac{1}{4r_j\ito\re^{2(u_j^0+w_j)}\ud x}\left[\tilde{Q}_j(s)-(-1)^{\vep_j}\sqrt{\tilde{Q}_j^2(s)-\frac{8b_j}{\lambda}\ito\re^{2(u_j^0+w_j)}\ud x}\right]
\nn\\
&&\equiv t_j-\varphi_{j,s}(t_{j-1},t_{j+1},\vep_j), \quad j=1, \dots, N,\label{a.3}
\eer
where
\ber
 \tilde{Q}_j(s)\equiv \ito\re^{u_j^0+w_j}(1+sr_{j-1}t_{j-1}\re^{u_{j-1}^0+w_{j-1}}+sr_{j+1}t_{j+1}\re^{u_{j+1}^0+w_{j+1}})\ud x.\label{a.3'}
\eer

We set,
\ber
 \Phi_{s,\vep}(t_1, \dots, t_N)=\big(\phi_{1,s}(t_1, t_2, \vep_1), \phi_{2,s}(t_1,t_2, t_3, \vep_2), \dots, \phi_{N,s}(t_{N-1},t_N, \vep_N)\big). \label{a.4}
\eer

 In what follows,   we always use $C$ to denote a universal positive constant whose value may change from line to line.
\begin{lemma}\label{lms21}
There exists a constant $C>1$, such that for given $\vep=(\vep_1, \dots, \vep_N)$ with $\vep_j\in\{0,1\},$ $s\in[0, 1]$,  and  $(t_1, \dots, t_N)$
 satisfying:  $ \Phi_{s,\vep}(t_1, \dots, t_N)=0$, we have:
 \ber
  \frac{1}{C}\frac{1}{\left(\ito\re^{2(u_j^0+w_j)}\ud x\right)^{1/2}}\le t_j\le   \frac{C}{\left(\ito\re^{2(u_j^0+w_j)}\ud x\right)^{1/2}},\, j=1, \dots, N. \label{1.4a}
 \eer
\end{lemma}

{\it Proof.} To establish \eqref{1.4a} we observe that, $t_j>0$, $\forall\,j=1,\dots, N$, and by setting:
 \[u_j=u_j^0+w_j+\ln t_j,\quad j=1, \dots, N,\]
we see that $u_j$ satisfies
\berr
&& 2r_j^2\ito\re^{2u_j}\ud x-sr_{j-1}r_j\ito\re^{u_{j-1}+u_j}\ud x-sr_jr_{j+1}\ito\re^{u_j+u_{j+1}}\ud x-r_j\ito\re^{u_j}\ud x\le 0, \nn\\
&&\quad \forall  j=1, \dots, N\,\text{and}\, s\in [0,1].
\eerr

 Since $M$  in \eqref{me} is positive definite,  with the help of H\"{o}lder's inequality and in view of the notation \eqref{a11}, we find  constants $\alpha_0>0, \beta_0>0$  such that,
\ber
  &&\alpha_0\sum\limits_{j=1}^N\ito\re^{2u_j}\ud x\le\ito\mathbf{U}^\tau M\mathbf{U}\ud x \le  \sum\limits_{j=1}^Nr_j\ito\re^{u_j}\ud x \nn\\
  &&\le |\Omega|^{\frac{1}{2}}\sum\limits_{j=1}^Nr_j\left(\ito\re^{2u_j}\ud x\right)^{\frac12}\le\beta_0\left(\sum\limits_{j=1}^N\ito\re^{2u_j}\ud x\right)^{\frac12}. \label{a.226}
\eer
Hence \eqref{a.226} implies that,
\ber
 \ito\re^{2u_j}\ud x\le C, \quad j=1, \dots, N \label{a.5}
\eer
 from which we readily get,
 \[
  t_j\le \frac{C}{\left(\ito\re^{2(u_j^0+w_j)}\ud x\right)^{1/2}},\quad j=1, \dots, N.
 \]

To obtain the reverse inequality, in view of \eqref{a.5}, we can estimate
\ber
\hat{Q}_j&=&\ito\re^{u_j^0+w_j}(1+sr_{j-1}t_{j-1}\re^{u_{j-1}^0+w_{j-1}}+sr_{j+1}t_{j+1}\re^{u_{j+1}^0+w_{j+1}})\ud x\nn\\
&\le& \left(\ito\re^{2u_j^0+2w_j}\ud x\right)^{\frac12}\left(|\Omega|^{\frac12}+r_{j-1}\left(\ito\re^{2u_{j-1}}\ud x\right)^{\frac12}+r_{j+1}\left(\ito\re^{2u_{j+1}}\ud x\right)^{\frac12}\right)\nn\\
&\le& C\left(\ito\re^{2u_j^0+2w_j}\ud x\right)^{\frac12},\label{a.6}
\eer
for suitable $C>0$ (depending only on $r_j, \, j=1, \dots, N$).

In case $\vep_j=1$,  then we can use \eqref{a.0} to derive:
\berr
t_j\ge \frac{\ito\re^{u_j^0+w_j}\ud x}{4r_j\ito\re^{2u_j^0+2w_j}\ud x}\ge \sqrt{\frac{8b_j}{\lambda}}\frac{1}{4r_j\left(\ito\re^{2u_j^0+2w_j}\ud x\right)^{\frac12}}
\eerr
and \eqref{a.4} is established in this case.

In case $\vep_j=0$,  then we can use \eqref{a.1} to deduce:
 \berr
   t_j&\ge& \frac{b_j}{\lambda r_j}\frac{1}{\ito\re^{u_j^0+w_j}(1+r_{j-1}\re^{u_{j-1}}+r_{j+1}\re^{u_{j+1}})\ud x}\\
   &\ge& \frac{b_j}{\lambda r_j}\frac{1}{C}\frac{1}{\left(\ito\re^{2u_j^0+2w_j}\ud x\right)^{\frac12}}
 \eerr
and \eqref{a.4} follows in this case as well.
\hfill $\square$\\

As a consequence of Lemma \ref{lms21}, we can take $R\gg1$ sufficiently large, such that the topological degree of $\Phi_{s,\vep}$
on $\Omega_R=\{(t_1, \dots, t_N): \, 0<t_j< R, \, j=1, \dots, N\}$ is well defined  for every $s\in [0, 1]$ and for every $\vep=(\vep_1, \dots, \vep_N)$ with
$\vep_j\in \{0, 1\}$, $j=1, \dots, N$.

By the homotopy invariance of the topological degree we find
\[
 \deg(\Phi_{s=1, \vep}, \Omega_R, 0)= \deg(\Phi_{s, \vep}, \Omega_R, 0)= \deg(\Phi_{s=0, \vep}, \Omega_R, 0)
\]
On the other hand, for any given $(w_1, \dots, w_N)\in \mathcal{A}$, we have:
$\Phi_{s=0, \vep}(t_1, \dots, t_N)=(t_1-a_1, \dots, t_N-a_N)$ with
\ber
 a_j&=&\frac{1}{4r_j\ito\re^{2(u_j^0+w_j)}\ud x}\times\nn\\
 &&\times\left(\ito\re^{u_j^0+w_j}\ud x+(-1)^{\vep_j}\sqrt{\left(\ito\re^{u_j^0+w_j}\ud x\right)^2-\frac{8b_j}{\lambda}\ito\re^{2(u_j^0+w_j)}\ud x}\right),\label{a.7}
\eer
$j=1, \dots, N$. Thus,  we obtain that,
\[
\deg(\Phi_{s=0, \vep}, \Omega_R, 0)=1.
\]

As a consequence,  $\deg(\Phi_{s=1, \vep}, \Omega_R, 0)= 1$ and we conclude that, for any given $(w_1, \dots, w_N)\in \mathcal{A}$  the system \eqref{a.1}
 admits at least one solution.

  To show that such a solution is actually unique and depends smoothly on  $(w_1, \dots, w_N)$, we show that $\forall\, s\in[0, 1]$ every solution of the
  equation:
   \be
   \Phi_{s, \vep}=0 \label{a.8}
   \ee
is actually non-degenerate. More precisely the following holds:
\begin{theorem}\label{thm21}
 For every  $(w_1, \dots, w_N)\in \mathcal{A}$, $\vep=(\vep_1, \dots, \vep_N)$, $\vep_j\in\{0,1\}, \, j=1, \dots, N$ and $s\in[0, 1]$, every solution
 of the equation \eqref{a.8} is nondegenerate.

\end{theorem}

{\it Proof.}  By the above calculation for $\Phi_{s=0,\vep}$, we see that the claim obviously holds for $s=0$. So we are left to consider the case where $0<s\le 1$.
To this purpose, we compute the Jacobian of $\Phi_{s,\vep}$. According to \eqref{a.3} and  \eqref{a.4}, we easily find that $\frac{\partial\phi_{j,s}}{\partial t_j}=1$
and
 \ber
\frac{\partial\Phi_{s,\vep}}{\partial t} =\begin{pmatrix}
 1&\frac{\partial\phi_{1,s}}{\partial t_2}& 0 &\cdots&\cdots&0\\
 \frac{\partial\phi_{2,s}}{\partial t_1}& 1& \frac{\partial\phi_{2,s}}{\partial t_3} &0& \cdots&0\\
  0&\frac{\partial\phi_{3,s}}{\partial t_2}& 1& \frac{\partial\phi_{3,s}}{\partial t_4} &\cdots&0\\
\vdots& &\ddots&\ddots&\ddots&\vdots\\
 \cdots&  & \ddots &\frac{\partial\phi_{N-1,s}}{\partial t_{N-2}}&1& \frac{\partial\phi_{N-1,s}}{\partial t_N}\\
 0& \cdots &\cdots&0&\frac{\partial\phi_{N,s}}{\partial t_{N-1}}&1
\end{pmatrix}.\label{a.4}
\eer

 As above, by setting $u_j=u_j^0+w_j+\ln t_j, \, j=1, \dots, N$, we find,
\ber
 &&\frac{\partial\phi_{j,s}}{\partial t_{j+1}}=-\frac{sr_{j+1}\ito\re^{u_j^0+w_j+u_{j+1}^0+w_{j+1}}\ud x}{4r_j\ito\re^{2(u_j^0+w_j)}\ud x}\times\nn\\
 &&\times\left[1-\frac{(-1)^{\vep_j}\ito\re^{u_j^0+w_j}(1+sr_{j-1}\re^{u_{j-1}}+sr_{j+1}\re^{u_{j+1}})\ud x}{\sqrt{\big(\ito\re^{u_j^0+w_j}(1+sr_{j-1}\re^{u_{j-1}}+sr_{j+1}\re^{u_{j+1}})\ud x\big)^2-\frac{8b_j}{\lambda}\ito\re^{2(u_j^0+w_j)}\ud x}}\right]. \label{a.10}
\eer

Therefore, if we use \eqref{a.3} we derive,
\ber
 \frac{\partial\phi_{j,s}}{\partial t_{j+1}}=\frac{(-1)^{\vep_j}st_jr_{j+1}\ito\re^{u_j^0+w_j+u_{j+1}^0+w_{j+1}}\ud x}{\sqrt{\big(\ito\re^{u_j^0+w_j}(1+sr_{j-1}\re^{u_{j-1}}+sr_{j+1}\re^{u_{j+1}})\ud x\big)^2-\frac{8b_j}{\lambda}\ito\re^{2(u_j^0+w_j)}\ud x}},\label{a.11}
\eer
$j=1, \dots, N$. Similarly, we find:
\ber
 \frac{\partial\phi_{j,s}}{\partial t_{j-1}}=\frac{(-1)^{\vep_j}st_jr_{j-1}\ito\re^{u_j^0+w_j+u_{j-1}^0+w_{j-1}}\ud x}{\sqrt{\big(\ito\re^{u_j^0+w_j}(1+sr_{j-1}\re^{u_{j-1}}+sr_{j+1}\re^{u_{j+1}})\ud x\big)^2-\frac{8b_j}{\lambda}\ito\re^{2(u_j^0+w_j)}\ud x}},\label{a.12}
\eer
$j=1, \dots, N$.

At this point, we are  going to use the (linear algebra) results of the Appendix A, in order to show that $ \det\frac{\partial \Phi_{s,\vep}}{\partial t}>0$.
To this purpose, from  \eqref{a.4}, \eqref{a.11} and \eqref{a.12}, we see that the Jacobian of $\Phi_{s,\vep}$ admits the same structure of the matrix
$T_1^{(N)}$ defined in \eqref{A.1} of the Appendix, with
\be
\beta_{j,1}=(-1)^{\vep_j}\alpha_{j,1}\quad \text{and}\quad \beta_{j,2}=(-1)^{\vep_j}\alpha_{j,2}, \quad j=1, \dots, N \quad \label{a.14}
\ee
and
\ber
 \alpha_{j,1}=\frac{st_{j-1}r_{j-1}\ito\re^{u_j^0+w_j+u_{j-1}^0+w_{j-1}}\ud x}{\sqrt{\big(\ito\re^{u_j^0+w_j}(1+sr_{j-1}\re^{u_{j-1}}+sr_{j+1}\re^{u_{j+1}})\ud x\big)^2-\frac{8b_j}{\lambda}\ito\re^{2(u_j^0+w_j)}\ud x}},\label{a.15a}\\
 \alpha_{j,2}=\frac{st_{j+1}r_{j+1}\ito\re^{u_j^0+w_j+u_{j+1}^0+w_{j+1}}\ud x}{\sqrt{\big(\ito\re^{u_j^0+w_j}(1+sr_{j-1}\re^{u_{j-1}}+sr_{j+1}\re^{u_{j+1}})\ud x\big)^2-\frac{8b_j}{\lambda}\ito\re^{2(u_j^0+w_j)}\ud x}}.\label{a.15b}
\eer


Therefore, the assumptions \eqref{A.7} and \eqref{A.8} of the Appendix are satisfied. Furthermore, concerning the coefficients $\alpha_{j,1}$ and $\alpha_{j,2}$, defined in \eqref{a.15a}--\eqref{a.15b} $j=1, \dots, N,$
we observe that, since $(w_1, \dots, w_N)\in \mathcal{A}$, then for $s\in(0,1]$ we can estimate:
\berr
 &&\left(\ito\re^{u_j^0+w_j}(1+st_{j-1}r_{j-1}\re^{u_{j-1}^0+w_{j-1}}+st_{j+1}r_{j+1}\re^{u_{j+1}^0+w_{j+1}})\ud x\right)^2\\
 &&>\frac{8b_j}{\lambda}\ito\re^{2(u_j^0+w_j)}\ud x+\left(st_{j-1}r_{j-1}\ito\re^{u_j^0+w_j+u_{j-1}^0+w_{j-1}}\ud x+st_{j+1}r_{j+1}\ito\re^{u_j^0+w_j+u_{j+1}^0+w_{j+1}}\ud x \right)^2.
\eerr

Consequently, by using the above estimate, for every s in (0, 1] we have:
\berr
0<\alpha_{j,2}<\frac{t_{j+1}r_{j+1}\ito\re^{u_j^0+w_j+u_{j+1}^0+w_{j+1}}\ud x}{t_{j-1}r_{j-1}\ito\re^{u_j^0+w_j+u_{j-1}^0+w_{j-1}}\ud x+t_{j+1}r_{j+1}\ito\re^{u_j^0+w_j+u_{j+1}^0+w_{j+1}}\ud x}\equiv1-\tau_j\le 1
\eerr
with
\berr
\tau_j=\frac{t_{j-1}r_{j-1}\ito\re^{u_j^0+w_j+u_{j-1}^0+w_{j-1}}\ud x}{t_{j-1}r_{j-1}\ito\re^{u_j^0+w_j+u_{j-1}^0+w_{j-1}}\ud x+t_{j+1}r_{j+1}\ito\re^{u_j^0+w_j+u_{j+1}^0+w_{j+1}}\ud x}\in[0, 1],
\eerr
and in turn,
\berr
0<\alpha_{j+1,1}<\frac{t_jr_j\ito\re^{u_j^0+w_j+u_{j+1}^0+w_{j+1}}\ud x}{t_jr_j\ito\re^{u_j^0+w_j+u_{j-1}^0+w_{j-1}}\ud x+t_{j+2}r_{j+2}\ito\re^{u_j^0+w_j+u_{j+2}^0+w_{j+2}}\ud x}\equiv \tau_{j+1}.
\eerr

In other words, for $j=1, \dots, N$, the coefficients  $\alpha_{j,1}$ and $\alpha_{j,2}$ in \eqref{a.15a}, \eqref{a.15b} satisfy also the assumption \eqref{A.18} of Theorem \ref{ath1} of the Appendix,  and  therefore we
may conclude that,
 \[\det\frac{\partial \Phi_{s,\vep}}{\partial t}=F^{(1)}_N>0\]
 and the  proof is completed.

\begin{corollary}\label{crls2}
 For every $(w_1, \dots, w_N)^\tau\in \mathcal{A}, \vep=(\vep_1, \dots, \vep_N)$ with $\vep_j\in\{0,1\}, j=1, \dots, N$ and $s\in[0, 1]$, the equation:
 \[
 \Phi_{s,\vep}(t_1, \dots, t_N)=0,
 \]
 admits a unique (non-degenerate) solution, smoothly depending on  $(w_1, \dots, w_N)$.
\end{corollary}

{\it Proof.}  It is clear that,  for $s=0$ the given statement follows form \eqref{a.7}. Furthermore, by the Implicit Function Theorem,  there exists $\delta>0$
sufficiently small,  such that for $s\in [0, \delta),$ problem \eqref{a.8} admits a unique (non-degenerate) solution. Let
\be
s_0\equiv\sup\left\{\sigma\in [0, 1]\Big|\, \text{ such that } \,\eqref{a.8} \, \text{ admits a unique solution }\quad \forall s\in [0, \sigma]\right\}.\label{a.s0}
\ee

We claim that,  $s_0=1$.

While it is clear that $s_0>0,$  if by contradiction, we suppose that $s_0<1$, then there would exist $s_0<s_n<1$ and $t_n^{(1)}\neq t_n^{(2)}\in \mathbb{R}^N$ such that,
\berr
&& \Phi_{s_n,\vep}(t_n^{(2)})= \Phi_{s_n,\vep}(t_n^{(1)})=0\\
&& s_n \searrow s_0  \quad \text{as}\quad n\to +\infty.
\eerr
 By virtue of Lemma \ref{lms21}, we can pass to a subsequence if necessary, to find that,
 \berr
  t_n^{(i)}\to t^{(i)} \quad \text{and} \quad \Phi_{s_0,\vep}(t^{(i)})=0,\quad i=1, 2.
 \eerr
By the  non-degeneracy of $t^{(i)}, i=1,2$, (as given by Theorem \ref{thm21}),  we can first rule out the possibility that, $t^{(1)}\neq t^{(2)}$. Indeed, if this was the case, then by the Implicit Function Theorem,
for sufficiently small $\delta>0$,  and for $s\in(0,1)$ such that: $s_0-\delta<s<s_0$ we would get that the equation \eqref{a.8} would admit at least two solutions, in contradiction with the definition of
$s_0$ in \eqref{a.s0}. Thus,   $t^{(1)}=t^{(2)}=\underline{t}$, and this would be again impossible, since the  Implicit Function Theorem implies local uniqueness
for solutions of \eqref{a.8} around $(s_0, \underline{t})$. \hfill $\square$

\section{Existence of multiple  solutions }\label{sct3}
\setcounter{equation}{0}\setcounter{lemma}{0}\setcounter{proposition}{0}\setcounter{remark}{0}\setcounter{theorem}{0}

 In this section we  show that  system \eqref{a14} admits at least two distinct solutions provided that  the parameter $\lambda$ is sufficiently  large and $N=3, 4, 5$.
 For this purpose,
by following  \cite{haya1,nota},  we  consider the  constrained functional
\ber
 J(\mathbf{w})\equiv I(\mathbf{w}+\mathbf{c}_+(\mathbf{w})),\quad  \mathbf{w}\in \mathcal{A}, \label{a40}
\eer
where $\mathbf{c}_+(\mathbf{w})$ is the unique solution of  the constraint equations \eqref{a26} with all  $\vep_j=1, \,  j=1,\dots, N,$ (see  Corollary \ref{crls2}).
By minimizing the the constrained functional $ J(\mathbf{w})$ in $\mathcal{A}$,  the authors of \cite{haya1} establish the following:
\begin{proposition}[\cite{haya1}]There exists $\overline{\lambda}>0$ such that for every $\lambda>\overline{\lambda}$ the functional $J$ in \eqref{a40} attains its minimum value at the point
$\mathbf{w}_\lambda$ which belongs to the $\bf {interior}$ of $ \mathcal{A}$. Namely,
 \[ J(\mathbf{w}_\lambda)=\inf\limits_{\mathcal{A}} J(\mathbf{w})\] and $\mathbf{v}^*_\lambda=  \mathbf{w}_\lambda + \mathbf{c}_+(\mathbf{w}_\lambda)$ defines a critical point for $I$  in  \eqref{a20}.
\end{proposition}
\hfill $\square$\\

By setting, $ \mathbf{c}_\lambda^*=\mathbf{c}_+(\mathbf{w}_\lambda)=(c_{1,\lambda}^*, \dots, c_{N,\lambda}^*)^\tau$
we may write,
\be
\mathbf{v}^*_\lambda=(v_{1,\lambda}^*, \dots, v_{N,\lambda}^*)^\tau=(c_{1,\lambda}^*+w_{1,\lambda}, \dots, c_{N,\lambda}^*+w_{N,\lambda})^\tau.\label{v32}
\ee
We will show that actually  $\mathbf{v}^*_\lambda$ is a local minimum for $I$ in  \eqref{a20}.
\begin{lemma}\label{lm3.1}
The solution $\mathbf{v}^*_\lambda=(v_{1,\lambda}^*, \dots, v_{N,\lambda}^*)^\tau$  of \eqref{a14}, as given by \eqref{v32},  defines a local minimum for functional $I$ in  \eqref{a20}.
\end{lemma}

 {\it Proof.}   For fixed $\mathbf{w}\in \mathcal{A}$, we denote by $\mathbf{c}^*(\mathbf{w})=(c_{1}^*(\mathbf{w}),\dots, c_{N}^*(\mathbf{w}))$  the unique solution of
 \eqref{a219} with $\vep_j=1, \,\forall j=1\dots, N,$ as given by Corollary \ref{crls2}.
    For $\mathbf{c}=(c_1, \dots, c_N)^\tau\in \mathbb{R}^N$ we easily check  that,
 \ber
\frac{\partial}{ \partial{c_j}}I(\mathbf{w}+\mathbf{c})\Big|_{\mathbf{c}=\mathbf{c}^*(\mathbf{w})}=\frac{\partial}{ \partial{c_j}}I(w_1+c_1, \dots, w_N+c_N)\Big|_{\mathbf{c}=\mathbf{c}^*(\mathbf{w})}=0, \quad j=1, \dots, N. \label{a3.3'}
 \eer

Moreover, by a straightforward computation we find:
 \ber
  &&\frac{\partial^2}{ \partial{c_j^2}}I(w_1+c_1, \dots, w_N+c_N)\nn\\
  &&=\lambda\ito r_j\re^{u_j^0+c_j+w_j}\big(4r_j\re^{u_j^0+c_j+w_j}-1-r_{j-1}\re^{u_{j-1}^0+c_{j-1}+w_{j-1}}-r_{j+1}\re^{u_{j+1}^0+c_{j+1}+w_{j+1}}\big)\ud x,\nn\\
  &&\quad j=1, \dots, N, \label{a3.4'}
 \eer
 and
 \ber
  &&\frac{\partial^2}{ \partial{c_jc_k}}I(w_1+c_1, \dots, w_N+c_N)=\frac{\partial^2}{ \partial{c_kc_j}}I(w_1+c_1, \dots, w_N+c_N),\nn\\
  &&=-\lambda r_jr_{k}\ito\re^{u_j^0+c_j+w_j+u_k^0+c_k+w_k}\ud x,\,\text{for }\, k\in\{ j-1, j+1\}\, \text{and }\, j=1, \dots, N,  \label{a3.5'}
 \eer
while,
\ber
  &&\frac{\partial^2}{ \partial{c_jc_k}}I(w_1+c_1, \dots, w_N+c_N)=\frac{\partial^2}{ \partial{c_kc_j}}I(w_1+c_1, \dots, w_N+c_N)=0,\nn\\
  &&\, \text{for}\,k\notin\{ j-1, j+1\}\, \text{and} \, j=1, \dots, N. \label{a3.6'}
 \eer

By setting $\mathbf{v}^*=\mathbf{w}+\mathbf{c}^*(\mathbf{w})=(v_1^*,\dots, v_N^*)^\tau$, then by the definition of $\mathbf{c}^*(\mathbf{w})$ and  \eqref{a219}, we see that,
\ber
&&\frac{\partial^2}{ \partial{c_j^2}}I(v_{1}^*, \dots, v_{N}^*)\nn\\
  &&=\lambda\left[\left(\ito r_j\re^{u_j^0+v^*_j}\big(1+r_{j-1}\re^{u_{j-1}^0+v^*_{j-1}}+r_{j+1}\re^{u_{j+1}^0+v^*_{j+1}}\big)\ud x\right)^2
  -\frac{8b_jr_j^2\ito\re^{2(u_j^0+v^*_j)}\ud x}{\lambda}\right]^{\frac12},\nn\\
  &&\quad \forall j=1, \dots, N.\label{a3.7'}
\eer

 In particular from \eqref{a3.4'}--\eqref{a3.7'} we  conclude  that,
\ber
 &&\frac{\partial^2}{ \partial{c_j^2}}I(v_{1}^*, \dots, v_{N}^*)+\frac{\partial^2}{ \partial{c_jc_{j-1}}}I(v_{1}^*, \dots, v_{N}^*)\nn\\
 &&+\frac{\partial^2}{ \partial{c_jc_{j+1}}}I(v_{1}^*, \dots, v_{N}^*)>0, \, \forall \,j=1, \dots, N.
 \label{a3.8'}
\eer
Hence from \eqref{a3.8'}  we infer  that, for any given $\mathbf{w}$  in  $\mathcal{A}$, the Hessian matrix of $I(\mathbf{w}+\mathbf{c}),$ as a function of $\mathbf{c}=(c_1, \dots, c_N),$  is a strictly diagonally dominant tri-diagonal matrix at the point $\mathbf{c}=\mathbf{c}^*(\mathbf{w}),$ and therefore it is  strictly positive-definite.
In particular this property holds for $\mathbf{w}_\lambda$ the minimun point of $J$ in \eqref{a40}.

Therefore  we conclude  that, for $\delta>0$ sufficiently small and  for $(v_1, \dots, v_N)^\tau=(w_1+c_1, \dots, w_N+c_N)^\tau$  satisfying:
  \be\sum\limits_{i=1}^N\|v_i-v^*_{i, \lambda}\|\le \delta, \label{3.9bis}\ee
 we have that  $(w_1, \dots,  w_N)^\tau$  belongs to  the interior of $\mathcal{A}$ and also
 (by the smooth dependence of $\mathbf{c}^*(\mathbf{w})$ with respect to $\mathbf{w}$ (see Corollary \ref{crls2})) that
 the vector $\mathbf{c}=(c_1, \dots, c_N)$ is sufficiently close to $\mathbf{c}^*(\mathbf{w})$ to guarantee that,
 \berr
   I(v_1, \dots, v_N)= I(w_1+c_1, \dots, w_N+c_N)&\ge& I(\mathbf{w}+\mathbf{c}^*(\mathbf{w}))=J(\mathbf{w})
  \eerr
As a consequence, for any $v$  satisfying \eqref{3.9bis}, we have:
\berr
   I(\mathbf{v})\ge \inf\limits_{\mathbf{w}\in \mathcal{A}}J(\mathbf{w})=I(\mathbf{v_\lambda}^*)
  \eerr
and the proof is completed.  \hfill $\square$\\

To proceed further, we need the following  ``compactness'' property of $I$.

 \begin{proposition} \label{props}
Let  $3\le N\le 5$ and  $\{(v_{1,n}, \dots, v_{N,n})^\tau\} \in (W^{1,2}(\Omega))^N$ be  such that,
\ber
 I(v_{1,n}, \dots, v_{N,n}) \to a_0\quad \text{as} \quad n\to +\infty,\label{a.3.10} \\
 \|I'(v_{1,n}, \dots, v_{N,n})\|_*\to 0 \quad \text{as} \quad n\to +\infty.,\label{a.3.11}
\eer
 where  $a_0$ is a constant and $\|\cdot\|_*$ denotes the norm of the dual space of  $(W^{1, 2}(\Omega))^N$. Then   $(v_{1,n}, \dots, v_{N,n})$  admits  a strongly convergent subsequence in $(W^{1, 2}(\Omega))^N$.

  \end{proposition}

Using a standard terminology, Proposition \ref{props} asserts that the functional $I$ satisfies the Palais--Smale (PS)-condition. We suspect that such property should hold also when $N\ge6.$

 We provide the proof of Proposition \ref{props} in the following section.

Based on Proposition \ref{props}, we can  carry out the proof of Theorem \ref{th1} and
 obtain a  second solution of  \eqref{a14} (other than  $(v^*_{1, \lambda}, \dots, v^*_{N, \lambda})^\tau$  in \eqref{v32})  by a Mountain-pass construction.

 To this purpose, we need to reduce to the case where we know that $\mathbf{v}^*_\lambda$  is a strict local minimum of $I$. Indeed, if on the contrary, for  small $\delta>0$,  we have:
\berr
 \inf\limits_{\sum\limits_{j=1}^N\|v_i-v^*_{i, \lambda}\|=\delta}I(v_1, \dots, v_N)=I(v^*_{1, \lambda}, \dots,  v^*_{N, \lambda}),
\eerr
then we conclude, from  Corollary 1.6 of \cite{gho}, that  the functional $I$ admits a  one parameter family of degenerate local minimizers, and
  a second solution of  \eqref{a14} is certainly guaranteed in  this case.

 Thus, we can assume that $\mathbf{v}^*_\lambda=(v^*_{1, \lambda}, \dots, v^*_{N, \lambda})^\tau$  is a  strict local minimum for $I$. So that  for   sufficiently small $\delta>0$, we have  that,
  \ber
   I(v^*_{1, \lambda}, \dots, v^*_{N, \lambda})<\inf\limits_{\sum\limits_{j=1}^N\|v_i-v^*_{i, \lambda}\|=\delta}I(v_1, \dots,  v_N)\equiv\gamma_0.\label{c.3.10}
  \eer
 On the other hand, we easily check that,
 \berr
  I(v^*_{1, \lambda}-\xi, \dots,  v^*_{N, \lambda}-\xi)\to -\infty\quad \text{as}\quad \xi\to +\infty.
 \eerr
  Therefore,   for a sufficiently large $\xi_0>1$,  we let
\berr
 \hat{v}_i\equiv v^*_{j, \lambda}-\xi_0, \quad j=1, \dots, N,
\eerr
and  conclude that,
 \ber
  \sum\limits_{j=1}^N\|\hat{v}_j-v^*_{j, \lambda}\|>\delta\label{c.3.11}
 \eer
 and
 \ber
  I(\hat{v}_1, \dots, \hat{v}_N)<I(v^*_{1, \lambda}, \dots, v^*_{N, \lambda})-1.\label{c.3.12}
 \eer

 We introduce the set of  paths,
\berr
\mathcal{P}\equiv\Big\{ \Gamma(t)\Big| \Gamma\in C\left([0, 1], \,  (W^{1, 2}(\Omega))^N\right), \quad \Gamma(0)=(v^*_{1, \lambda}, \dots,  v^*_{N, \lambda})^\tau, \quad
\Gamma(1)=(\hat{v}_1, \dots, \hat{v}_2)^\tau\Big\}
\eerr
and define:
 \berr
  a_0\equiv\inf\limits_{\Gamma\in\mathcal{P}}\sup\limits_{t\in[0, 1]}I(\Gamma(t)).\label{c.3.14}
 \eerr
Clearly,
 \be
 a_0> I(v^*_{1, \lambda}, \dots, v^*_{N, \lambda}).\label{c.3.13}
  \ee
and in view of Proposition \ref{props}, we can use  the ``Mountain-pass'' theorem of  Ambrosetti-Rabinowitz  \cite{amra} to obtain that
 $a_0$ defines  a critical  value of the functional $I$, to which it corresponds to  a  critical point  different from $(v^*_{1, \lambda}, \dots,  v^*_{N, \lambda})^\tau$.
  Thus the proof of Theorem \ref{th1} is completed. \hfill $\square$\\


\section{The (PS)-condition for $N=3, 4, 5.$} \label{sct4}
\setcounter{equation}{0}\setcounter{lemma}{0}\setcounter{proposition}{0}\setcounter{remark}{0}\setcounter{theorem}{0}

We devote this section to establish the (PS)-condition.

Let $\{(v_{1,n}, \dots, v_{N,n})\}$ be a sequence in $(W^{1,2}(\Omega))^N$ satisfying \eqref{a.3.10}--\eqref{a.3.11} and denote by,
\[u_{j,n}=u_j^0+v_{j,n},\quad j=1, \dots, N, \]
where $u_j^0$ is given by \eqref{a.2.11}.  In what follows,  we always use the decomposition:
\be
v_{j,n}=c_{j,n}+w_{j,n}, \quad w_{j,n}\in \dot{W}^{1, 2}(\Omega), \, c_{j,n}=\frac{1}{|\Omega|}\ito v_{j,n}\ud x,\, j=1, \dots, N, \forall\,n\in \mathbb{N}. \label{a4.0}
\ee

By recalling  \eqref{re} and \eqref{re2} , we note that,
\be
 2r_j-r_{j-1}-r_{j+1}=1\quad\forall\, j=1, \dots, N, \label{4.*}
\ee
and \eqref{4.*}, allow us to obtain the following:
 \ber
 && I'(v_{1,n},\dots, v_{N,n})(\phi_1,\dots, \phi_N)\nn\\
 &&= \sum\limits_{j,k=1}^Na_{jk}\ito\nabla w_{k,n}\cdot\nabla \phi_j\ud x+\lambda\sum\limits_{j=1}^N\ito r_j\re^{u_{j,n}}\Big(2r_j(\re^{u_{j,n}}-1)-r_{j-1}(\re^{u_{j-1,n}}-1)\nn\\
  &&\quad-r_{j+1}(\re^{u_{j+1,n}}-1)\Big)\phi_j\ud x+\frac{1}{|\Omega|}\sum\limits_{j=1}^N\ito b_j\phi_j\ud x\nn\\
&&= \sum\limits_{j,k=1}^Na_{jk}\ito\nabla w_{k,n}\cdot\nabla \phi_j\ud x+\lambda\sum\limits_{j=1}^N\ito r_j\re^{u_{j,n}}(2r_j\re^{u_{j,n}}-r_{j-1}\re^{u_{j-1,n}}-r_{j+1}\re^{u_{j+1,n}})\phi_j\ud x\nn\\
  &&\quad -\lambda\sum\limits_{j=1}^Nr_j\ito \re^{u_{j,n}}\phi_j\ud x+\frac{1}{|\Omega|}\sum\limits_{j=1}^N\ito b_j\phi_j\ud x\nn\\
&&=o(1)\|\phi\|, \quad \forall \mathbf{\phi}=(\phi_1, \dots, \phi_N)\in (W^{1,2}(\Omega))^N.\label{1.1}
 \eer

Thus, by taking $(\phi_1,\dots, \phi_N)=(1,0, \dots, 0), \dots, (0,\dots, 0, 1)$ separately in \eqref{1.1}, we find:
\ber
 &&\lambda\ito\Big[2r_j^2\re^{2u_{j,n}}-r_jr_{j-1}\re^{u_{j,n}+u_{j-1,n}}-r_jr_{j+1}\re^{u_{j,n}+u_{j+1,n}}-r_j\re^{u_{j,n}}\Big]\ud x+b_j\nn\\
 &&=\lambda\ito\Big[2r_j^2\re^{u_{j,n}}(\re^{u_{j,n}}-1)-r_jr_{j-1}\re^{u_{j,n}}(\re^{u_{j-1,n}}-1)-r_jr_{j+1}\re^{u_{j,n}}(\re^{u_{j+1,n}}-1)\Big]\ud x+b_j\nn\\
 &&=o(1),\, j=1, \dots, N.\label{a4.5}
\eer
Hence by using still \eqref{4.*}, we can sum up the identities \eqref{a4.5} over $j=1, \dots, N,$ and arrive at the following identity:
\ber
&&\lambda\ito(\mathbf{U}_n-\mathbf{1})^\tau M(\mathbf{U}_n-\mathbf{1}) \ud x+\lambda\sum\limits_{j=1}^Nr_j\ito\re^{u_{j,n}}\ud x\nn\\
&&-\frac{N(N+1)(N+2)}{12}\lambda|\Omega|+\sum\limits_{j=1}^Nb_j= o(1),\label{a4.6}
\eer
with \[\mathbf{U}_n=(\re^{u_{1,n}}, \dots, \re^{u_{N,n}})^\tau.\]
Since  the matrix $M$ is positive definite,  from \eqref{a4.6} we see that, as $n\to +\infty$:
\ber
 &&\ito(\mathbf{U}_n-\mathbf{1})^\tau M(\mathbf{U}_n-\mathbf{1}) \ud x\le \frac{N(N+1)(N+2)}{12}|\Omega|+o(1),\label{a4.7}\\
 &&\ito\re^{u_{j,n}}\ud x\le  \frac{N(N+1)(N+2)}{12}|\Omega|+o(1),\quad j=1, \dots, N. \label{a4.8}
\eer

By Jensen's inequality, from \eqref{a4.0}, \eqref{a4.7} and \eqref{a4.8},  we have:
\ber
 \re^{c_{j,n}}\le \frac{N(N+1)(N+2)}{12}+o(1), \quad j=1, \dots, N. \label{a4.9}
\eer

In addition, from \eqref{a4.7} and \eqref{a4.8},   we derive that
  \ber
   \int_\Omega(\re^{u_j^0+v_{j, n}}-1)^2\ud x\le C, \quad j=1, \dots, N, \label{a4.10}\\
   \int_\Omega\re^{2u_j^0+2v_{j, n}} \ud x\le C, \quad j=1, \dots, N.\label{a4.11}
  \eer
 for some suitable constant $C>0$.

Therefore, if we take $(\phi_1,\dots, \phi_N)=(w_{1,n}, \dots, w_{N,n})$ in \eqref{1.1},  in view of \eqref{a4.9}--\eqref{a4.11},  we find positive constants $\beta_1>0$ and $\beta_2>0$   such that,
 \ber
  &&o(1)\sum\limits_{i=1}^N\|\nabla w_{j,n}\|_2\nm\\
  &\ge&I'(v_{1,n}, \dots, v_{N,n})[(w_{1,n}, \dots, w_{N,n})]\nm\\
  &=& \sum\limits_{j,k=1}^N a_{ij}\ito\nabla w_{j,n}\cdot\nabla w_{j,n}\ud x+\lm\sum\limits_{j=1}^N\ito r_j\re^{u_{j,n}}(2r_j\re^{u_{j,n}}-r_{j-1}\re^{u_{j-1,n}}-r_{j+1}\re^{u_{j+1,n}})w_{j,n}\ud x\nn\\
  &&\quad -\lambda\sum\limits_{j=1}^Nr_j\ito \re^{u_{j,n}}w_{j,n}\ud x\nm\\
  &\ge& \beta_1 \sum\limits_{j=1}^N\|\nabla w_{j,n}\|_2^2+\lm\sum\limits_{j=1}^N\ito 2r_j^2\re^{u_{j,n}}\left[\re^{u_j^0+c_{j,n}}(\re^{w_{j,n}}-1)+\re^{u_j^0+c_{j,n}}\right]w_{j,n}\ud x\nm\\
  &&-\lm\sum\limits_{j=1}^N\ito r_jr_{j+1}\re^{u_{j,n}+u_{j+1,n}}(w_{j,n}+w_{j+1,n})\ud x-\lambda\sum\limits_{j=1}^Nr_j\ito \re^{u_{j,n}}w_{j,n}\ud x\nm\\
  &\ge& \beta_1 \sum\limits_{j=1}^N\|\nabla w_{j,n}\|_2^2-\beta_2\sum\limits_{j=1}^N\|\nabla w_{j,n}\|_2-\lm\sum\limits_{j=1}^N r_jr_{j+1}\ito\re^{u_{j,n}+u_{j+1,n}}(w_{j,n}+w_{j+1,n})\ud x.\label{f1'}
 \eer

 At this point, our main effort  will be to obtain an uniform  estimate for the term:\\
 $\ito\re^{u_{j,n}+u_{j+1,n}}(w_{j,n}+w_{j+1,n})\ud x,$  \, for every $ j=1, \dots, N$.

We start by showing that $w_{j,n}$ is uniformly bounded in $L^p,$ for any $p > 1$.

To this purpose,  for fixed $j\in\{1,\cdots,N\}$  we take:
 \[\phi_k=\left\{\begin{array}{lll}
  0&k\notin\{j-1, j, j+1\}\\
  -\varphi&k=j-1, j+1\\
  2\varphi&k=j
  \end{array}\right.
,\]
and from \eqref{1.1} we obtain:
\ber
&&\ito \nabla w_{j,n}\cdot \nabla\varphi \ud x+\lambda\ito\Big[4r_j^2\re^{2u_{j,n}}-2r_{j-1}^2\re^{2u_{j-1,n}}-2r_{j+1}^2\re^{2u_{j+1,n}}+r_{j-1}r_{j-2}\re^{u_{j-1,n}+u_{j-2,n}}
\nn\\
&&+r_{j+1}r_{j+2}\re^{u_{j+1,n}+u_{j+2,n}}-r_{j-1}r_j\re^{u_{j-1,n}+u_{j,n}}-r_jr_{j+1}\re^{u_{j,n}+u_{j+1,n}}+2r_j\re^{u_{j, n}}-r_{j-1}\re^{u_{j-1,n}}\nn\\
&&-r_{j+1}\re^{u_{j+1,n}}\Big]\varphi\ud x+\frac{1}{|\Omega|}(2b_j-b_{j-1}-b_{j+1})\ito\varphi\ud x=o(1)\|\varphi\|,\quad \forall \varphi\in W^{1,2}(\Omega). \label{1.3}
\eer

For any  $1<q<2$,  we know that,
\ber
\|\nabla w_{j,n}\|_q=\sup \left\{\ito \nabla w_{j,n}\cdot \nabla\varphi \ud x,  \quad  \forall\varphi\in W^{1, p}(\Omega): \, \ito\varphi\ud x =0, \|\nabla\varphi\|_p=1;  \frac1p+\frac1q=1\right\}\label{1.4}
\eer
and every $\varphi$ in \eqref{1.4} satisfies: $\|\varphi\|+\|\varphi\|_\infty\le C$, for suitable $C > 0.$ \\
Thus, from \eqref{1.3} and \eqref{1.4}   we derive:
\[ \|\nabla w_{j,n}\|_q\le C_q, \, \text{for some }\, C_q>0. \]

As a consequence, for any $p\ge 1$, there exists $C_p>0$ such that:
\ber
 \|w_{j,n}\|_p\le C_p, \quad  \forall n\in \mathbb{N}\,\text{ and } \quad j=1,\dots,N. \label{1.5}
\eer

Next,  in \eqref{1.1} we take $\phi_j=\varphi\in W^{1,2}(\Omega)$ for every $ \, j=1, \dots, N$, and  by simple calculations, we get,
\berr
 &&\ito\nabla\Big(\sum\limits_{j,k=1}^Na_{jk}w_{j,n}\Big)\cdot\nabla \varphi\ud x +\lambda \sum\limits_{j=0}^N\ito(r_j\re^{u_{j,n}}-r_{j+1}\re^{u_{j+1,n}})^2\varphi\ud x
\nn\\
&& +\lambda\sum\limits_{j=1}^N\left(-r_j\ito\re^{u_{j,n}}\varphi\ud x+\frac{b_j}{|\Omega|}\ito \varphi\ud x\right)= o(1)\|\varphi\|,
\eerr
and since  $\sum\limits_{k=1}^Na_{jk}=r_j$, we find:
\ber
&& \ito\nabla\Big(\sum\limits_{j=1}^Nr_jw_{j,n}\Big) \cdot\nabla\varphi\ud x+ \lambda\ito\left[\sum\limits_{j=0}^N(r_j\re^{u_{j,n}}-r_{j+1}\re^{u_{j+1,n}})^2-\sum\limits_{j=0}^Nr_j\re^{u_{j,n}}\right]\varphi\ud x
 \nn\\
 &&+\sum\limits_{j=1}^N\frac{b_j}{|\Omega|}\ito \varphi\ud x =  o(1)\|\varphi\|,\quad \forall \varphi\in W^{1,2}(\Omega).\label{1.2}
 \eer

Therefore, we can choose  $\varphi=\big(\sum\limits_{j=1}^Nr_jw_{j,n}\big)^+$ in \eqref{1.2},  and   in view of \eqref{1.5} we find,
\ber
 \Big\|\nabla\Big(\sum\limits_{j=1}^Nr_jw_{j,n}\Big)^+\Big\|_2\le C.\label{1.6a}
\eer

As usual, we have denoted by $f^+(x)=\max\{f(x),0\}$ the positive part of $f=f(x)$.

More generally we define:
\[
 W_n^{(j)}=\sum\limits_{k=j}^Nr_kw_{k,n}\quad\text{and}\quad  \hat{W}_n^{(j)}=\sum\limits_{k=1}^{N+1-j}r_kw_{k,n},
\]
with $r_j$ in  \eqref{re}--\eqref{re2}, and show the following:
\begin{lemma} \label{lm5.1}
 If $3\le N\le 5$, then
 \ber
  \|\nabla(W_n^{(j)})^+\|_2\le C\quad\text{and}\quad  \|\nabla(\hat{W}_n^{(j)})^+\|_2\le C, \quad \forall j=1, \dots, N. \label{2.3}
 \eer
\end{lemma}

{\it Proof.}  If $j=1$ then \eqref{2.3} reduces to \eqref{1.6a}. Hence we take $j\ge 2$ and proceed by induction. Note that, by the symmetry:
\ber
 r_j=r_{N+1-j},\quad j\in\left\{1, \dots, \Big[\frac{N+1}{2}\Big]\right\}. \label{2.0}
\eer
it suffices to prove the uniform estimate for $W_n^{(j)}$ as the one for $\hat{W}_n^{(j)}$ follows similarly. \\
Next, we observe that,
 \be
 \forall\, j=1,\dots, N,\quad  (2+r_j)^2 < 8(1+r_j) \quad \text{if and only if} \quad 1\le N\le 5,  \label{2.1}
 \ee
and \eqref{2.1} is the exact reason for which we need the  restriction on $N\in \{3, 4, 5\}.$

To check \eqref{2.1}, observe that it is equivalent to:
\ber
 r_j<2(1+\sqrt2)\quad \forall j=1,\dots, N,\label{2.1'}
\eer

and by  \eqref{2.0} it suffices to check it only for $1\le j\le [\frac{N+1}{2}]$. But for such $j$'s the value of $r_j$
is increasing with respect to $j$,  and so \eqref{2.1} holds if and only if $r_{[\frac{N+1}{2}]}< 2(1+\sqrt2).$

If $N=2k$ is even  then $[\frac{N+1}{2}]=k$ and $r_k=\frac{k(k+1)}{2},$ while for $N=2k+1$  odd we have: $[\frac{N+1}{2}]=k+1$ and $r_{k+1}=\frac{(k+1)^2}{2}$.
Hence \eqref{2.1'} holds if and only if $k=0,1,2$, namely, $N=1,2,3,4,5$, as claimed.

Next, we illustrate the induction scheme for  $j=2$, where we  use \eqref{1.1} with $\phi_1=-r_2\varphi$, $\phi_2=(1+r_1)\varphi$ and $\phi_k=\varphi \,\, \forall\, k=3,\dots, N$. We obtain:
\ber
 &&\ito\nabla W_n^{(2)}\cdot\nabla\varphi\ud x+\lambda\ito\left(-2r_2r_1^2\re^{2u_{1,n}}+(r_2-(1+r_1))\re^{u_{1,n}+u_{2,n}}\right)\varphi\ud x\nn\\
 &&+\lambda\ito(2(1+r_1)r_2^2\re^{2u_{2,n}}-(2+r_1)r_2r_3\re^{u_{2,n}+u_{3,n}}+r_3^2\re^{2u_{3,n}})\varphi\ud x\nn\\
 &&+\lambda\ito\left(\sum\limits_{k=3}^N(r_k\re^{u_{k,n}}-r_{k+1}\re^{u_{k+1,n}})^2\right)\varphi\ud x=o(1)\|\varphi\|.\label{2.4}
\eer

Thus, by using \eqref{2.1} with $j=1$, we let:
 \[\mu_1=\left(2(1+r_1)-\frac{(2+r_1)^2}{4}\right)r_2^2>0\]
and rewrite \eqref{2.4} as follows:
\berr
&&\ito \nabla W_n^{(2)}\cdot\nabla\varphi\ud x+\lambda\ito\left(-2r_2r_1^2+\frac{(r_2-(1+r_1))^2}{4\mu_1}\right)\re^{2u_{1,n}}\varphi\ud x\nn\\
&&+\lambda\ito\left[\left(\frac{r_2-(1+r_1)}{2\sqrt{\mu_1}}\re^{u_{1,n}}+\sqrt{\mu_1}\re^{u_{2,n}}\right)^2+\left(\frac{2+r_1}{2}r_2\re^{u_{2,n}}-r_3\re^{u_{3,n}}\right)^2\right]\varphi\ud x\\
 &&+\lambda\ito\left(\sum\limits_{k=3}^N(r_k\re^{u_{k,n}}-r_{k+1}\re^{u_{k+1,n}})^2\right)\varphi\ud x=o(1)\|\varphi\|.
\eerr

Therefore, by choosing $\varphi=(W_n^{(2)})^+$ in \eqref{2.4} and by recalling \eqref{a4.9},
 we obtain:
\berr
&&\|\nabla(W_n^{(2)})^+\|_2^2\le C\left(\ito\re^{2w_{1,n}}(W_n^{(2)})^+\ud x+1\right)\nn\\
&&\le C\left(\int_{\Omega\cap\{w_{1,n}\ge0\}}\re^{\frac{2}{r_1}(r_1w_{1,n}+(W_n^{(2)})^+)}(W_n^{(2)})^+\ud x+1\right)\nn\\
&&\le  C\left(\int_{\Omega\cap\{w_{1,n}\ge0\}}\re^{\frac{2}{r_1}(\sum\limits_{j=1}^Nr_jw_{j,n})^+}(W_n^{(2)})^+\ud x+1\right)\nn\\
&&\le  C\left(\ito\re^{\frac{2}{r_1}(\sum\limits_{j=1}^Nr_jw_{j,n})^+}(W_n^{(2)})^+\ud x+1\right).
\eerr

At this point, by virtue of \eqref{1.5}  and  \eqref{1.6a}, we can use  the Moser--Trudinger inequality
(see \cite{aubi}), to conclude that,
 \[\ito\re^{\frac{2}{r_1}(W_n^{(1)})^+}(W_n^{(2)})^+\ud x\le C,\]
and   we arrive at the desired conclusion.

Next, let us carry out our induction procedure, so  for $j\ge 3$ we assume that,
\ber
 \|\nabla(W_n^{(j-1)})^+\|_2\le C,\label{2.4'}
\eer
 and we are left to  prove that a similar estimate holds for $(W_n^{(j)})^+$.

To this purpose, we are going to use \eqref{1.1} with $\phi_{j-1}=-r_j\varphi, \phi_j=(1+r_{j-1})\varphi, \phi_k=0$ for $1\le k\le j-2$ and $\phi_k=\varphi$ for
$k\ge j+1,$ and obtain the following:
\berr
 &&\ito\nabla W_n^{(j)}\cdot\nabla\varphi \ud x+\lambda\ito\big(-2r_jr^2_{j-1}\re^{2u_{j-1,n}}+(r_j-(1+r_{j-1}))\re^{u_{j-1,n}+u_{j,n}}\big)\varphi\ud x\\
 &&+\lambda\ito\big(2(1+r_{j-1})r_j^2\re^{2u_{j,n}}-(2+r_{j-1})r_jr_{j+1}\re^{u_{j,n}+u_{j+1,n}}+r_{j+1}^2\re^{2u_{j+1,n}}\big)\varphi\ud x\\
 &&+\lambda\ito\left(\sum\limits_{k=j+1}^N(r_k\re^{u_{k,n}}-r_{k+1}\re^{u_{k+1,n}})^2\right)\varphi\ud x=o(1)\|\varphi\|.
\eerr

Exactly as above, by using \eqref{2.1} and by letting:
\be
 \mu_{j-1}=\left(2(1+r_{j-1})-\frac{(2+r_{j-1})^2}{4}\right)r_2^2>0, \quad \sigma_j=\frac{(r_j-(1+r_{j-1}))^2}{4\mu_{j-1}}+2r_jr^2_{j-1}>0,\label{2.5}
\ee
we find:
\berr
&&\ito \nabla W_n^{(j)}\cdot\nabla\varphi \ud x+\lambda\ito\left(-\sigma_j\re^{2u_{j-1,n}}
+\left(\frac{r_j-(1+r_{j-1})}{2\sqrt{\mu_{j-1}}}\re^{u_{j-1,n}}+\sqrt{\mu_{j-1}}\re^{u_{j,n}}\right)^2\right)\varphi\ud x\\
&&+\lambda\ito\left(\frac{2+r_{j-1}}{2}r_j\re^{u_{j,n}}-r_{j+1}\re^{u_{j+1,n}}\right)^2\varphi\ud x\\
&&+\lambda\ito\left(\sum\limits_{k=j+1}^N\left(r_k\re^{u_{k,n}}-r_{k+1}\re^{u_{k+1,n}}\right)^2\right)\varphi\ud x=o(1)\|\varphi\|.
\eerr

Therefore, by taking $\varphi=(W_n^{(j)})^+$ we derive:
\berr
&&\|\nabla(W_n^{(j)})^+\|_2^2\le C\left(\ito\re^{2w_{j-1,n}}(W_n^{(j)})^+\ud x+1\right)\nn\\
&&\le C\left(\int_{\Omega\cap\{w_{j-1,n}\ge0\}}\re^{\frac{2}{r_{j-1}}(r_{j-1}w_{j-1,n}+(W_n^{(j)})^+)}(W_n^{(j)})^+\ud x+1\right)\nn\\
&&\le  C\left(\ito\re^{\frac{2}{r_{j-1}}(W_n^{(j-1)})^+}(W_n^{(j)})^+\ud x+1\right)\le C,
\eerr
and the last estimate follows as above,  from the induction hypothesis \eqref{2.4}, \eqref{1.5} and the Moser--Trudinger inequality \cite{aubi}.

Again by the symmetry \eqref{2.0}, we can argue similarly simply by replacing  in the above arguments each index involved, say $k$, with $N+1-k$ and obtain:
 \[
  \|\nabla (\hat{W}_n^{(j)})^+\|_2\le C,\quad \forall j=1, \dots, N
 \]
and the proof is completed.

\hfill $\square$\\

At this point, we observe that:  $r_1w_{1,n}+r_2w_{2,n}=\hat{W}_n^{(N-1)}$ and $r_1w_{1,n}=\hat{W}_n^{(N)}$;  while,
 $r_{N-1}w_{N-1,n}+r_Nw_{N,n}=W_n^{(N-1)}$ and $r_Nw_{N,n}=W_n^{(N)}$. Thus,  from Lemma \ref{lm5.1}, we obtain  in particular  that,
 \ber
  \|\nabla(r_1w_{1,n}+r_2w_{2,n})^+\|_2+\|\nabla(r_{N-1}w_{N-1,n}+r_Nw_{N,n})^+\|_2\le C, \label{2.5a}
 \eer
 and
 \ber
  \|\nabla w_{1,n}^+\|_2+\|\nabla w_N^+\|_2\le C. \label{2.5b}
 \eer

 More generally there holds:
 \begin{lemma} \label{lmf.2}
  For $N=3, 4, 5$ we have:
  \be
   \|(r_jw_{j,n}+r_{j+1}w_{j+1,n})^+\|\le C, \quad \forall\, j=1, \dots, N. \label{5.*}
  \ee
 \end{lemma}

{\it Proof.}  For $N=3$, \eqref{5.*} follows already from \eqref{1.5}, \eqref{2.5a} and \eqref{2.5b}. So we let,
 $N\ge4$ and set,
\be
 Z_n=\sum\limits_{k=2}^{N-1}r_kw_{k,n}.\label{2.6}
\ee

{\bf CLAIM 1:}

 \be
 \|\nabla(Z_n)^+\|_2\le C.\label{2.7}
 \ee

  To establish \eqref{2.7}  we apply \eqref{1.1} with $\phi_1=-r_2\varphi$, $\phi_2=(1+r_1)\varphi$, $\phi_{N-2}=(1+r_{N-3})\varphi$, $\phi_{N-1}=(1+r_N)\varphi$ and
$\phi_N=-r_{N-1}\varphi$. Note that for $N=4$ the definition above is consistent since $\phi_2=(1+r_1)\varphi=(1+r_{N-3})\varphi=\phi_{N-2}$.

For $N=4$  we have that, $r_1=r_4=2$ and $r_2=r_3=3,$ so  we obtain:
\berr
&&\ito\nabla Z_n\cdot\nabla\varphi \ud x+\lambda\ito\Big[-2r_1^2r_2\re^{2u_{1,n}}-(1+r_1-r_2)r_1r_2\re^{u_{1,n}+u_{2,n}}+2(1+r_1)(r_2\re^{u_{2,n}}-r_3\re^{u_{3,n}})^2
\nn\\&&-2r_3r_4^2\re^{2u_{4,n}}+(1+r_4-r_3)r_3^2r_4\re^{u_{3,n}+u_{4,n}}\Big]\varphi\ud x\nn\\
&&=\ito\nabla Z_n\cdot\nabla\varphi\ud x-24\lambda\ito(\re^{2u_{1,n}}+\re^{2u_{4,n}})\varphi\ud x+6\lambda\ito(r_2\re^{u_{2,n}}-r_3\re^{u_{3,n}})^2\varphi\ud x\\
&&=o(1)\|\varphi\|.
\eerr

Hence, by taking $\varphi=(Z_n)^+$, and by using \eqref{2.5b} and the Moser--Trudinger inequality as above,  we derive:
\berr
\|\nabla (Z_n)^+\|_2^2\le C\ito(\re^{2w_{1,n}}+\re^{2w_{4,n}})(Z_n)^+\ud x\le
C\ito(\re^{2w_{1,n}^+}+\re^{2w_{4,n}^+})(Z_n)^+\ud x\le C
\eerr
 and \eqref{2.7} follows  for $N=4$.

For $N=5$,  we have:
\be
r_1=r_5=\frac52, \, r_2=r_4=4, \, r_3=\frac92, \label{2.8}
\ee
and in this case from \eqref{1.1} we obtain:
\berr
&&\ito\nabla Z_n\cdot\nabla\varphi\ud x+\lambda\ito\big[-2r_1^2r_2\re^{2u_{1,n}}-(1+r_1-r_2)r_1r_2\re^{u_{1,n}+u_{2,n}}\big]\varphi\ud x
\\&&+\lambda\ito\left[2(1+r_1)r_2^2\re^{2u_{2,n}}-(2+r_1+r_2)r_2r_3\re^{u_{2,n}+u_{3,n}}+2(1+r_2)r_3^2\re^{2u_{3,n}}\right]\varphi\ud x
 \\&&+\lambda\ito\big[2(1+r_5)r_4^2\re^{2u_{4,n}}-(2+r_2+r_5)r_3r_4\re^{u_{3,n}+u_{4,n}}-2r_4r_5^2\re^{2u_{5,n}}\\
 &&+r_4r_5(r_4-(1+r_5))\re^{u_{4,n}+u_{5,n}}\big]\varphi\ud x\\
&&=\ito\nabla Z_n\cdot\nabla\varphi\ud x-18\lambda\ito(\re^{2u_{1,n}}+\re^{2u_{5,n}})\varphi\ud x
+5\lambda\ito(\re^{u_{1,n}+u_{2,n}}+\re^{u_{4,n}+u_{5,n}})\varphi\ud x\\
&&+\lambda\ito[2(1+r_1)r_2^2\re^{2u_{2,n}}-(2+r_1+r_2)r_2r_3\re^{u_{2,n}+u_{3,n}}+(1+r_2)r_3^2\re^{2u_{3,n}}]\varphi\ud x
 \\&&+\lambda\ito[2(1+r_1)r_2^2\re^{2u_{4,n}}-(2+r_1+r_2)r_2r_3\re^{u_{3,n}+u_{4,n}}+(1+r_2)r_3^2\re^{2u_{3,n}}]\varphi\ud x\\
&&=o(1)\|\varphi\|.
\eerr

At this point, we observe that,  $(2+r_1+r_2)^2\le 8(1+r_1)(1+r_2)$, and so   we can  check that the terms within the brackets in the last two integrals above are positive.

 As a consequence  for $\varphi=(Z_n)^+,$ arguing as above we find:
\berr
\|\nabla (Z_n)^+\|_2^2\le C\ito(\re^{2w_{1,n}^+}+\re^{2w_{5,n}^+})(Z_n)^+\ud x\le C,
\eerr
 and  Claim 1 is established for $N=5$ as well.

 {\bf CLAIM 2:} If $N=5$ then,
\ber
 \|\nabla(r_jw_{j,n}+r_{j+1}w_{j+1,n})^+\|_2\le C,\quad j= 2, 3.\label{*}
\eer

To establish \eqref{*}, we observe that,    from Lemma \ref{lm5.1} and  \eqref{2.7}, we have:
\ber
 \|\nabla(r_jw_{j,n}+r_{j+1}w_{j+1,n}+r_{j+2}w_{j+2,n})^+\|_2\le C,\quad j=1, 2, 3.\label{2.9}
\eer
So, for $j=2,3$, we can take $\phi_{j-1}=-r_j\varphi$, $\phi_j=(1+r_{j-1})\varphi$, $\phi_{j+1}=(1+r_{j+2})\varphi$,  $\phi_{j+2}=-r_{j+1}\varphi$,
in \eqref{1.1} and obtain:
\ber
 &&\ito \nabla(r_jw_{j,n}+r_{j+1}w_{j+1,n})\cdot\nabla\varphi\ud x+\lambda\ito\Big[-2r_{j-1}^2r_j\re^{2u_{j-1,n}}\nn\\
 &&-(1+r_{j-1}-r_j)r_{j-1}r_j\re^{u_{j-1,n}+u_{j,n}}+2(1+r_{j-1})r_j^2\re^{2u_{j,n}}-(2+r_{j-1}+r_{j+1})r_jr_{j+1}\re^{u_{j,n}+u_{j+1,n}}\nn\\
 &&+2(1+r_{j+2})r_{j+1}^2\re^{2u_{j+1,n}}-(1+r_{j+2}-r_{j+1})r_{j+1}r_{j+2}\re^{u_{j+1,n}+u_{j+2,n}}\nn\\
 &&-2r_{j+1}r_{j+2}^2\re^{2u_{j+2,n}}+r_{j+1}r_{j+2}r_{j+3}\re^{u_{j+2,n}+u_{j+3,n}}\Big]\varphi\ud x=o(1)\|\varphi\|.\label{2.10}
\eer

Now notice that, for $j=2, 3$ we have:
\[1+r_{j-1}-r_j=\frac{2j-N}{2}=-(1+r_{j+2}-r_{j+1})\]
and by \eqref{2.8} we can check directly that,
\ber
 (2+r_{j-1}+r_{j+2})^2<16(1+r_{j-1})(1+r_{j+2}). \label{2.11}
\eer

Hence \eqref{2.10} can be expressed as follows:
\berr
&&\ito\nabla(r_jw_{j,n}+r_{j+1}w_{j+1,n})\cdot \nabla\varphi\ud x-2\lambda\ito(r_{j-1}^2r_j\re^{2u_{j-1,n}}+r_{j+1}r_{j+2}^2\re^{2u_{j+2,n}})\varphi\ud x\\
 &&+\lambda\frac{N-2j}{2}\ito(r_{j-1}r_j\re^{u_{j-1,n}+u_{j,n}}-r_{j+1}r_{j+2}\re^{u_{j+1,n}+u_{j+2,n}})\varphi\ud x\\
 &&+\lambda\ito\Big[2(1+r_{j-1})r_j^2\re^{2u_{j,n}}-(2+r_{j-1}+r_{j+2})r_jr_{j+1}\re^{u_{j,n}+u_{j+1,n}}+2(1+r_{j+2})r_{j+1}^2\re^{2u_{j+1,n}}\Big]\varphi\ud x\\
 &&+\lambda r_{j+1}r_{j+2}r_{j+3}\ito\re^{u_{j+2,n}+u_{j+3,n}}\varphi\ud x=o(1)\|\varphi\|.
\eerr

Consequently, for $j=2$, we set:
\[
\vep_1=\Big(2(1+r_4)-\frac{(2+r_1+r_4)^2}{8(1+r_1)}\Big)r_3^2>0
\]
(recall \eqref{2.11}) and find:
\berr
 &&\ito \nabla(r_2w_{2,n}+r_3w_{3,n})\cdot\nabla\varphi\ud x-2\lambda\ito\Big(r_1^2r_2\re^{2u_{1,n}}+r_3r_4^2\Big(1+\frac{r_3}{32\vep_1}\Big)\re^{2u_{4,n}}\Big)\varphi\ud x\\
  &&+\lambda\ito\Big(\frac12r_1r_2\re^{u_{1,n}+u_{2,n}}+r_3r_4r_5\re^{u_{4,n}+u_{5,n}}\Big)\varphi\ud x\\
  &&+\lambda\ito\Big(\sqrt{2(1+r_1)}r_2\re^{u_{2,n}}-\frac{2+r_1+r_4}{2\sqrt{2(1+r_1)}}r_3\re^{u_{3,n}}\Big)^2\varphi\ud x\\
  &&+\lambda\ito\Big(\sqrt\vep_1\re^{u_{3,n}}-\frac{r_3r_4}{4\sqrt\vep_1}\re^{u_{4,n}}\Big)^2\varphi\ud x=o(1)\|\varphi\|.
\eerr

Thus, by taking $\varphi=(r_2w_{2,n}+r_3w_{3,n})^+$, we can estimate:
\berr
&&\|\nabla(r_2w_{2,n}+r_3w_{3,n})^+\|_2^2\\
&&\le C\left(\ito\re^{2w_{1,n}}(r_2w_{2,n}+r_3w_{3,n})^+\ud x +\ito\re^{2w_{4,n}}(r_2w_{2,n}+r_3w_{3,n})^+\ud x\right)+C_1\\
&&\le C\left(\int_{\{w_{1,n}\ge0\}}\re^{r_1w_{1,n}}(r_2w_{2,n}+r_3w_{3,n})^+\ud x +\int_{\{w_{4,n}\ge0\}}\re^{r_4w_{4,n}}(r_2w_{2,n}+r_3w_{3,n})^+\ud x\right)+C_2\\
&&\le C\left(\ito\big(\re^{(r_1w_{1,n}+r_2w_{2,n}+r_3w_{3,n})^+}+\re^{(r_2w_{2,n}+r_3w_{3,n}+r_4w_{4,n})^+}\big)(r_2w_{2,n}+r_3w_{3,n})^+\ud x\right)+C_2\\
 &&\le C
\eerr
as it follows from \eqref{2.9} and again by using  \eqref{1.5} together with  the Moser--Trudinger inequality. Thus \eqref{*} is established for $j=2$.

For $j=3$, we argue similarly and for,
\[
\vep_2=\Big(2(1+r_2)-\frac{(2+r_2+r_5)^2}{8(1+r_5)}\Big)r_3^2>0
\]
we obtain:
 \berr
 &&\ito \nabla(r_3w_{3,n}+r_4w_{4,n})\cdot\nabla\varphi\ud x-2\lambda\ito\Big[r_2^2r_3\Big(1+\frac{r_3}{32\vep_2}\Big)\re^{2u_{2,n}}+r_4r_5^2\re^{2u_{5,n}}\Big]\varphi\ud x\\
 &&+\frac\lambda2r_4r_5\ito\re^{u_{4,n}+u_{5,n}}\varphi\ud x+\lambda\ito\Big(r_3\frac{2+r_2+r_5}{2\sqrt{2(1+r_5)}}\re^{u_{3,n}}-\sqrt{2(1+r_5)}r_4\re^{u_{4,n}}\Big)^2\varphi\ud x\\
 &&+\lambda\ito(\sqrt{\vep_2}\re^{u_{3,n}}-\frac{r_2r_3}{4\sqrt{\vep_2}}\re^{u_{2,n}})^2\varphi\ud x=o(1)\|\varphi\|.
 \eerr
As above,  for $\varphi=(r_3w_{3,n}+r_4w_{4,n})^+$ we  obtain the following estimate:
\berr
 &&\|\nabla(r_3w_{3,n}+r_4w_{4,n})^+\|_2^2\le C\ito(\re^{2w_{2,n}}+\re^{2w_{5,n}})(r_3w_{3,n}+r_4w_{4,n})^+\ud x+C_1\\
&&\le C\left(\int_{\{w_{2,n}\ge0\}}\re^{r_2w_{2,n}}(r_3w_{3,n}+r_4w_{4,n})^+\ud x +\int_{\{w_{5,n}\ge0\}}\re^{r_5w_{5,n}}(r_3w_{3,n}+r_4w_{4,n})^+\ud x\right)+C_2\\
&&\le C\left(\ito(\re^{(r_2w_{2,n}+r_3w_{3,n}+r_4w_{4,n})^+}+\re^{(r_3w_{3,n}+r_4w_{4,n}+r_5w_{5,n})^+})(r_3w_{3,n}+r_4w_{4,n})^+\ud x\right)+C_2\\
 &&\le C,
\eerr
and \eqref{*} is established.  Thus, the proof of \eqref{5.*} is completed.
 \hfill $\square$\\
\\

 {\it The  Proof of Proposition \ref{props}: }
\\
\\
 According to  \eqref{a4.9} and Lemma \eqref{lmf.2}, for $j=1, \dots, N$, we can estimate:
 \berr
&&\ito\re^{u_{j,n}+u_{j+1,n}}(w_{j,n}+w_{j+1,n})\ud x\le C_1 \, (\int\limits_{\{w_{j,n}\ge0\}\}\cap\{w_{j+1,n}\ge0\}}\re^{w_{j,n}+w_{j+1,n}}(w_{j,n}+w_{j+1,n})\ud x  \, \,)+C_2\nn\\
&&\le C_1 \,(\int\limits_{\{w_{j,n}\ge0\}\}\cap\{w_{j+1,n}\ge0\}}\re^{(r_jw_{j,n}+r_{j+1}w_{j+1,n})^+}(w_{j,n}+w_{j+1,n})\ud x \, \, )+C_2\le C,
\eerr
where the last estimates follows as above, by \eqref{1.5}, Lemma \ref{lmf.2} and the Moser--Trudinger inequality.

 At this point, by virtue of \eqref{f1'} we may conclude also  that,
 \[ \sum\limits_{j=1}^N\|\nabla w_{j,n}\|_2^2\le C.\]

 In particular, along a subsequence, as $n\to +\infty$, we may conclude that,
 \berr
  w_{j,n}\to w_j\quad \text {weakly in }\, W^{1,2}(\Omega) \, \text{ and strongly in  } \, L^p(\Omega),\\
  \re^{w_{j,n}}\to \re^{ w_j} \quad \text{ strongly in  }\quad L^p(\Omega), \quad \text{for }\, p>1.
 \eerr
 Moreover, from \eqref{a4.5}, we  derive that, (along a subsequence) there holds:
 \[  \re^{2c_{j,n}}\ito\re^{2u_j^0+2w_{j,n}}\ud x\to L_j>0, \quad \text{as}\quad n\to +\infty,\]
with a suitable $ L_j > 0. $ \\
Consequently, by setting:
$$c_j = \frac{1} {2} \, log ( \frac{L_j}{\ito\re^{2u_j^0+2w_{j}}\ud x} ) \text{ and }  v_j=w_j+c_j, $$
$\text { as }  n\to +\infty, $ we have:


 \[ v_{j,n}\to v_j, \, \quad  \re^{u_j^0+v_{j,n}} \to \re^{u_j^0+v_j}\, \, \text{ strongly in  }\,L^p(\Omega) \,\text{ for }\, p>1, \,\quad \forall j=1, \dots, N.\]

Therefore, from \eqref{1.2},  we get:
 \berr
  &&\sum\limits_{j,k=1}^N a_{jk}\ito \nabla(w_{j,n}-w_j)\cdot\nabla(w_{k,n}-w_k)\ud x \\
 && =I'(v_{1,n}, \dots, v_{N,N})(v_{1,n}-v_1, \dots, v_{N,n}-v_N)+o(1)\\
 &&=o(1)\left(\sum\limits_{j=1}^N \|\nabla(w_{j,n}-w_j)\|_2+1\right), \quad \text{as}\, n\to+\infty
 \eerr

which implies that,
\berr
  w_{j,n}\to w_j\quad \text{strongly in}\quad W^{1,2}(\Omega), \quad \text{as}\, n\to +\infty, \quad  \forall j=1, \dots,N,\\
 \eerr
and the proof is completed.

\appendix\renewcommand{\appendixname}{Appendix~\Alph{section}}
 \section{Appendix of Linear Algebra}
\setcounter{equation}{0}\setcounter{lemma}{0}\setcounter{proposition}{0}\setcounter{remark}{0}\setcounter{theorem}{0}

For  $N\ge 2$ and for $j=1, \dots, N$ we let,
\ber
 \beta_{j,1}\in \mathbb{R} \quad \text{for}\,j=2, \dots, N, \quad \beta_{j,2}\in \mathbb{R} \quad \text{for}\, j=1, \dots, N-1 \, \text{and set}\, \beta_{1,1}=0=\beta_{N,2}.\label{A.0}
\eer
Given $1\le k\le l\le N$, we define the square $(l-k+1)\times(l-k+1)$ matrix $T_l^{(k)}$ as follows:
 \ber
T_l^{(k)}=(t_{j,s})_{j,s=k,\dots, l}\quad \text{with}\quad t_{j,s}=\beta_{j,1}\delta_{j-1}^s+\delta_j^s+\beta_{j,2}\delta_{j+1}^s \label{A.1}
 \eer
where,
 \[ \delta_p^s=\left\{\begin{array}{lll}
 1, &s=p\\
 0, &s\neq p
 \end{array} \right.
  \]
 is the usual Kronecker symbol.

Notice that, if $k=l$ then $T_l^{(l)}=1,$  while in general the matrix $T_l^{(k)}$ is expressed in terms of the quantities:
\ber
 (\beta_{k,2}, \beta_{k+1,1},  \beta_{k+1,2}, \dots,  \beta_{l-1,1},\beta_{l-1,2}, \beta_{l,1}). \label{A.2}
\eer

We wish to identify the determinant of $T^{(k)}_l,$  i.e.  $\det T^{(k)}_l.$ \\
To this purpose  we define the quantities $F_l^{(k)}$ via a recursive formula, starting with the case $k=l$, where we set,
\be
 F_l^{(l)}=1. \label{A.3}
\ee

For $k>l$  we define $F_l^{(k)}$ as follows:
\be
 F_l^{(l+1)}=1\quad \text{and} \quad F_l^{(k)}=0\quad \forall\, k \ge l+2. \label{A.4}
\ee

More importantly, for $1\le k\le l$ recursively we set:
\ber
 &&F_l^{(k)}=1-\sum\limits_{j=k}^{l}\beta_{j+1,1}\beta_{j,2}F_l^{(j+2)}\nn\\
 &&= 1-\sum\limits_{j=k}^{l}\beta_{j+1,1}\beta_{j,2}\left(1-\sum\limits_{k=j+2}^{l}\beta_{k+1,1}\beta_{k,2}\left(1-\sum\limits_{s=k+2}^{l}\beta_{s+1,1}\beta_{s,2}(1-\cdots)\right)\right)\label{A.5}
\eer
so that,
\[
 F_l^{(l-1)}=1-\beta_{l,1}\beta_{l-1,2},\quad  F_l^{(l-2)}=1-\beta_{l-2,2}\beta_{l-1,1}-\beta_{l-1,2}\beta_{l,1}
\]
and so on.

\begin{remark}
In view of \eqref{A.4},  it suffices to take the summation in \eqref{A.5} up to the index  $k = l-1$, instead of $l$ as indicated there.
\end{remark}

Notice that the value of  $F_l^{(k)}$ depends  on the same terms of the matrix $T_l^{(k)},$ as specified in \eqref{A.2}.  In fact, the following holds:
 \begin{lemma}\label{Al1}
 For $1\le k\le l\le N$ we have:
 \be
  \det T_l^{(k)}=F_l^{(k)}. \label{A.6}
 \ee
 \end{lemma}

 {\it Proof.}  We proceed by induction on $k$. In fact for $k=l,$ in view of \eqref{A.4}, we see that \eqref{A.6} is obviously satisfied.
 Thus, we assume that for $1\le k\le l-1,$ there holds:
 \berr
  \det T_l^{(j)}=F_l^{(j)}\quad \forall\, j\in \{k+1, \dots, l\}.
 \eerr
To check \eqref{A.6}, we observe that,
\[
 \det T_l^{(k)}= \det T_l^{(k+1)}-\beta_{k+1,1}\beta_{k,2} \det T_l^{(k+2)}.
\]
Thus,  by the induction hypothesis,    we find that,
\[
 \det T_l^{(k)}= F_l^{(k+1)}-\beta_{k+1,1}\beta_{k,2}F_l^{(k+2)}=F_l^{(k)}
\]
and the proof is complete.
\hfill $\square$
\\

To proceed further, we need to be more specific about our choice of $\beta$'s in \eqref{A.0}.   More precisely we let,
\ber
 &&\beta_{j,1}=(-1)^{\vep_{j,1}}\alpha_{j,1}, \, \vep_{j,1}\in\{0,1\}, \alpha_{j,1}\ge0,\, j=2, \dots, N,\label{A.7}\\
 &&\beta_{j,2}=(-1)^{\vep_{j,2}}\alpha_{j,2}, \, \vep_{j,2}\in\{0,1\}, \alpha_{j,2}\ge0,\, j=1, \dots, N-1, \, \text{and}\, \alpha_{1,1}=0=\alpha_{N,2}.\label{A.8}
\eer

\begin{lemma}\label{Al2}
Assume \eqref{A.7} and \eqref{A.8}. Then, for given $1\le k\le l\le N$  there holds:
\ber
 \frac{\partial F_l^{(k)}}{\partial\alpha_{j+1,1}}=-(-1)^{\vep_{j,2}+\vep_{j+1,2}}\alpha_{j,2}F_l^{(j+2)}F_{j-1}^{(k)}, \label{A.9}\\
  \frac{\partial F_l^{(k)}}{\partial\alpha_{j,2}}=-(-1)^{\vep_{j,2}+\vep_{j+1,1}}\alpha_{j+1,1}F_l^{(j+2)}F_{j-1}^{(k)}.\label{A.10}
\eer
\end{lemma}

{\it Proof.}  First of all,  by virtue of \eqref{A.4}, we can check that, for $1\le j<k$ we have: $F_{j-1}^{(k)}=0.$ Similarly, for $l\le j\le N$ we have $F_l^{j+2}=0.$
 So, for such choice of indices, we have:  $ \frac{\partial F_l^{(k)}}{\partial\alpha_{j+1,1}}=0=\frac{\partial F_l^{(k)}}{\partial\alpha_{j,2}}$
consistently with the definition of $F_l^{(k)}$.  Hence we let, $1\le k<l$ and for $j\in \{k, \dots, l-1\}$ we are going to verify \eqref{A.9}  and \eqref{A.10}
 by an induction argument on $k$. Actually, we provide the details only for \eqref{A.9}, as \eqref{A.10} follows similarly.

 For $k=l-1$ we see that,
 \[
  F_{l-1}^{(l-1)}=1-(-1)^{\vep_{l-1,1}+\vep_{l,1}}\alpha_{l-1,2}\alpha_{l,1},
 \]
and in this case  we have only the choice of $j=k=l-1$.   Hence,
\[
\frac{\partial F_l^{(k)}}{\partial\alpha_{j+1,1}}=\frac{\partial F_{l-1}^{(l-1)}}{\partial\alpha_{l,1}}=(-1)^{\vep_{l-1,1}+\vep_{l,1}}\alpha_{l-1,2},
\]
which gives exactly \eqref{A.9}, since in this case we  have:
\[
F_l^{(j+2)}=F_l^{(l+2)}=1 \quad \text{and}\quad F_{j-1}^{(k)}=F_{l-2}^{(l-1)}=1.
\]

Next, we take $k\in\{1, \dots, l-2\}$ and  by induction we assume that, for $j\in\{k+1, \dots, l-1\}$ the identity \eqref{A.9} holds for $\frac{\partial F_l^{(k)}}{\partial\alpha_{j+1,1}}.$

For $k\le j<l,$ we write,
\berr
&&F_l^{(k)}=1-\sum\limits_{s=k}^{j-1}(-1)^{\vep_{s,2}+\vep_{s+1,1}}\alpha_{s,2}\alpha_{s+1,1}F_l^{(s+2)}-(-1)^{\vep_{j,2}+\vep_{j+1,1}}\alpha_{j,2}\alpha_{j+1,1}F_l^{(j+2)}\\
&&-\sum\limits_{s=j+1}^{l}(-1)^{\vep_{s,2}+\vep_{s+1,1}}\alpha_{s,2}\alpha_{s+1,1}F_l^{(s+2)}
\eerr
with the understanding that, when  $j=k$ then the first summation term above is dropped. We compute:
\berr
 \frac{\partial F_l^{(k)}}{\partial\alpha_{j+1,1}}=-\sum\limits_{s=k}^{j-1}(-1)^{\vep_{s,2}+\vep_{s+1,1}}\alpha_{s,2}\alpha_{s+1,1}\frac{\partial F_l^{(s+2)}}{\partial\alpha_{j+1,1}}
 -(-1)^{\vep_{j,2}+\vep_{j+1,1}}\alpha_{j,2}F_l^{(j+2)}.
\eerr

  Hence, by virtue of the induction assumption we find:
  \berr
  &&\frac{\partial F_l^{(k)}}{\partial\alpha_{j+1,1}}=-\left(1-\sum\limits_{s=k}^{j-1}(-1)^{\vep_{s,2}+\vep_{s+1,1}}\alpha_{s,2}\alpha_{s+1,1}F_{j-1}^{(s+2)}\right)(-1)^{\vep_{j,2}+\vep_{j+1,1}}\alpha_{j,2}F_l^{(j+2)}\\
  &&=-(-1)^{\vep_{j,2}+\vep_{j+1,1}}\alpha_{j,2}F_l^{(j+2)}F_{j-1}^{(k)}
  \eerr
as claimed.
\hfill $\square$

\begin{remark}\label{ark1}
  Note that the term $F_l^{(j+2)}F_{j-1}^{(k)}$ on the right-hand side of \eqref{A.9} and \eqref{A.10} is independent of $\alpha_{j,2}$ and $\alpha_{j+1,1}$. Therefore if such
  term vanishes then $F_l^{(k)}$ is independent of both  $\alpha_{j,2}$ and $\alpha_{j+1,1}$.
\end{remark}

\begin{proposition}\label{ap1}
Let $i\in\{1, \dots, N-1\}$ and let $k,l\in \mathbb{N}$ be such that $k\le i\le l \le N$. We have:
\ber
&& \text{i) if}\quad  \vep_{i,2}+\vep_{i+1,1}=0\,\text{ ( mod2 ) }\quad \text{ then}\quad   F_l^{(k)}=F_l^{(k)}|_{\vep_{i,2}=0=\vep_{i+1,1}};\label{A.11}\\
&& \text{ii) if}\quad \vep_{i,2}+\vep_{i+1,1}=1 \text{ and}\quad  F_l^{(i+2)}F_{i-1}^{(k)}\ge0\quad \text{then}\quad  F_l^{(k)}\ge F_l^{(k)}|_{\vep_{i,2}=0=\vep_{i+1,1}}.\label{A.12}
\eer
\end{proposition}

Although it is intuitively clear, we wish to clarify the notation adopted in \eqref{A.11} and \eqref{A.12} before presenting the proof of  Proposition \ref{ap1}.
We have set,
\be
 F_l^{(k)}\big|_{\vep_{i,2}=0=\vep_{i+1,1}}=\det(\overline{T}_{l,i}^{(k)}),\label{A.13}
\ee
where,
\ber
&&\overline{T}_{l,i}^{(k)}= (\overline{t}_{j,s}^i)_{j,s=k,\dots,l} \quad \text{with} \label{A.14}\\
&&\overline{t}_{j,s}^i=t_{j,s}\,\text{(defined in \eqref{A.1}) for}\quad  j\notin\{i, i+1\},  \label{A.15}\\
&&\overline{t}_{i,s}^i=(-1)^{\vep_{i,1}}\alpha_{i,1}\delta_{i-1}^s+\delta_i^s+\alpha_{i,2}\delta_{i+1}^s, \label{A.16}\\
&&\overline{t}_{i+1,s}^i=\alpha_{i+1,1}\delta_i^s+\delta_{i+1}^s+(-1)^{\vep_{i+1,2}}\alpha_{i+1,2}\delta_{i+2}^s.\label{A.17}
\eer
\\

{\it Proof of Proposition \ref{ap1} }
\\
\\
In order to establish \eqref{A.11}, we proceed by induction on $k$. Indeed, if $k=i$ then
\[
  F_l^{(i)}= F_l^{(i+1)}-(-1)^{\vep_{i,2}+\vep_{i+1,1}}\alpha_{i,2}\alpha_{i+1,1}F_l^{(i+2)}=F_l^{(i+1)}-\alpha_{i,2}\alpha_{i+1,1}F_l^{(i+2)},
\]
as follows  by the assumption:  $\vep_{i,2}+\vep_{i+1,1}=0$ ( mod2 ).

Since neither $F_l^{(i+1)}$ nor $F_l^{(i+2)}$ depends on the terms $\alpha_{i,2}$ and $\alpha_{i+1,1}$ (recall \eqref{A.3}), we see that,
\[F_l^{(i)}=F_l^{(i)}|_{\vep_{i,2}=0=\vep_{i+1,1}}.\]

Next,  for $1\le k< i$, suppose that
 \[F_l^{(s)}=F_l^{(s)}|_{\vep_{i,2}=0=\vep_{i+1,1}}, \quad \forall\,s\in \{k+1, \dots, l\}.\]

To establish that the same identity also holds for $s=k$, we observe that,
\berr
&&F_l^{(k)}=F_l^{(k+1)}-(-1)^{\vep_{k,2}+\vep_{k+1,1}}\alpha_{k,2}\alpha_{k+1,1}F_l^{(k+2)}\\
&&=F_l^{(k+1)}\big|_{\vep_{i,2}=0=\vep_{i+1,1}}-\alpha_{k,2}\alpha_{k+1,1}F_l^{(k+2)}\big|_{\vep_{i,2}=0=\vep_{i+1,1}}=F_l^{(k)}\big|_{\vep_{i,2}=0=\vep_{i+1,1}},
\eerr
where the last identity is a consequence of the definition in \eqref{A.13}--\eqref{A.17}.

To prove \eqref{A.12}, we use the derivation formulae \eqref{A.9} and \eqref{A.10} which under the given assumptions imply that $F_l^{(k)}$
is increasing separately with respect to $\alpha_{i,2}$ and $\alpha_{i+1,1}$.  In other words, in case $\vep_{i,2}=0$ and $\vep_{i+1}=1$ we have that:
\berr
&& F_l^{(k)}((-1)^{\vep_{k,2}}\alpha_{k,2}, \dots, \alpha_{i,2}, -t, (-1)^{\vep_{i+1,2}}\alpha_{i+1,2}, \dots, (-1)^{\vep_{l,1}}\alpha_{l,1})\\
 &&\ge F_l^{(k)}((-1)^{\vep_{k,2}}\alpha_{k,2}, \dots, \alpha_{i,2}, t, (-1)^{\vep_{i+1,2}}\alpha_{i+1,2}, \dots, (-1)^{\vep_{l,1}}\alpha_{l,1})
\eerr
and \eqref{A.12} follows by taking $t=\alpha_{j+1,1}$.  Similarly, if $\vep_{i,2}=1$ and $\vep_{i+1}=0$ then,
\berr
&& F_l^{(k)}((-1)^{\vep_{k,2}}\alpha_{k,2}, \dots, \alpha_{i,1}, -t, \alpha_{i+1,1}, \dots, \alpha_{l,1})\\
 &&\ge F_l^{(k)}((-1)^{\vep_{k,2}}\alpha_{k,2}, \dots, \alpha_{i,1}, t, \alpha_{i+1,1}, \dots, \alpha_{l,1}), \quad \forall\, t\in\mathbb{R},
\eerr
and in this case \eqref{A.12} follows by taking $t=\alpha_{i,2}$.

\hfill $\square$
\\
\\
The main purpose of this Appendix is  to establish the following result:
\begin{theorem}\label{ath1}
 Let $1\le k\le l\le N$ and assume that \eqref{A.7} and \eqref{A.8} hold. \\
 For given $\tau_j\in[0, 1]  \quad j=1, \dots, l$, we suppose that,
 \be
  0\le \alpha_{j,2}<1-\tau_j\quad \text{and}\quad 0\le \alpha_{j+1,1}<\tau_{j+1}\quad j=k, \dots, l-1, \label{A.18}
 \ee
 then
 \be
  \det T_l^{(k)}=F_l^{(k)}>0. \label{A.19}
 \ee
\end{theorem}

The proof will be given in several steps. Firstly,  we proceed  to prove \eqref{A.19} in case:
\be
 \vep_{j,2}=0=\vep_{j+1,1}\quad \forall, j=k, \dots, l-1. \label{A.20}
\ee
Thus, we let $T_{l,0}^{(k)}$ be the matrix defined in \eqref{A.1}--\eqref{A.2} with $\vep_{j,i}$ satisfying \eqref{A.20}, namely
 \ber
T_{l,0}^{(k)}=\begin{pmatrix}
 1&\alpha_{k,2}& 0 &0&\cdots&0\\
 \alpha_{k+1,1}& 1& \alpha_{k+1,2}&0& \cdots&0\\
 0& \alpha_{k+2,1}&1&\alpha_{k+2,2}&0&0\\
 \vdots&  & \ddots &\ddots &\ddots&\vdots\\
 0& \quad  & \ddots &\alpha_{l-1,1}&1&\alpha_{l-1,2}\\
 0& \cdots&\cdots   &0&\alpha_{l,1}&1
\end{pmatrix},\label{A.21}
\eer
and set
\be
F_{l,0}^{(k)}=\det T_{l,0}^{(k)}.\label{A.22}
\ee

Furthermore, for $0\le \tau_j\le 1,$ and $j=k, \dots, l$,  we introduce the matrix:
\be
\overline{T}_l^{(k)}=(\overline{t}_{js})_{j,s=k,\dots, l} \label{A.23}
\ee
with
\be
\overline{t}_{js}=\tau_j\delta_{j-1}^s+\delta_j^s+(1-\tau_j)\delta_{j+1}^s\label{A.24}
\ee
and set,
\be
\overline{F}_l^{(k)}=\det\overline{T}_l^{(k)}=1-\sum\limits_{j=k}^l\tau_{j+1}(1-\tau_j)\overline{F}_l^{(j+2)}. \label{A.25}
\ee
Concerning $\overline{F}_l^{(k)}$  we have the following:
\begin{lemma}\label{al3}
 Let $1\le k\le l\le N$ then
 \be
 \overline{F}_l^{(k)}(\tau_k, \dots, \tau_l)\ge 0, \quad \forall\,\tau_j\in[0, 1],\, j=k, \dots, l.\label{A.26}
 \ee
 Moreover,
  \be
 \overline{F}_l^{(k)}(0,\tau_{k+1}, \dots, \tau_{l-1}, 1)=0.\label{A.27}
 \ee
\end{lemma}

{\it Proof. } We can establish \eqref{A.27} simply by observing that $\overline{F}_l^{(k)}(0,\tau_{k+1}, \dots, \tau_{l-1}, 1)$ coincides with the determinant of
the matrix,
\[
\begin{pmatrix}
 1&1& 0 &0&\cdots&0\\
 \tau_{k+1}& 1& 1-\tau_{k+1}&0& \cdots&0\\
 0&\tau_{k+2}& 1& 1-\tau_{k+2}&0& 0\\
 \vdots&  & \ddots &\ddots &\ddots&\vdots\\
 0& & \ddots &\tau_{l-1}&1&1-\tau_{l-1}\\
 0& \cdots & &0&1&0
\end{pmatrix},
\]
which is clearly singular. Indeed, the sum of the odd columns coincides with the sum of the even columns and it is given by the column with all
entries equal to 1. Hence we have the linear dependence of the column-vectors and \eqref{A.27} follows.

To establish \eqref{A.26} we proceed by induction on $n\in \mathbb{N}$, with $0\le l-k\le n$. Indeed, for $n=1$ then either $k=l$ and $\overline{F}_l^{(l)}=1$
or $k=l-1$ and
\[\overline{F}_l^{(l-1)}(\tau_{l-1}, \tau_l)=1-\tau_l(1-\tau_{l-1})\ge 0.\]

Thus, we assume that $n>1$ and suppose that,
\be
\overline{F}_m^{(s)}(\tau_m, \dots, \tau_s)\ge 0\quad \forall\, m,s\in \mathbb{N}: \,  0\le m-s\le n-1.\label{A.29}
\ee
Hence, for $l\ge 3$  and $l-k=n,$ we need to prove that  $\overline{F}_l^{(k)}\ge 0$.  To this purpose we observe that,
\berr
\overline{F}_l^{(k)}=\overline{F}_l^{(k+1)}-\tau_{k+1}(1-\tau_k)\overline{F}_l^{(k+2)}=
\overline{F}_{l-1}^{(k+1)}-\tau_l(1-\tau_{l-1})\overline{F}_{l-2}^{(k+1)}-\tau_{k+1}(1-\tau_k)\overline{F}_l^{(k+2)}.
\eerr
According to our induction assumption \eqref{A.29}, we see that $\overline{F}_{l-2}^{(k+1)}\ge0$ and $\overline{F}_l^{(k+2)}\ge0$ and therefore,
\[
\overline{F}_l^{(k)}\ge \overline{F}_{l-1}^{(k+1)}-\tau_{k+1}\overline{F}_l^{(k+2)}-(1-\tau_{l-1})\overline{F}_{l-2}^{(k+1)}=
\overline{F}_l^{(k)}(0, \tau_{k+1}, \dots, \tau_{l-1}, 1)=0
\]
and \eqref{A.26} is established.

\hfill $\square$\\

\begin{lemma} \label{al4}
Let $1\le k<l\le N$ and assume \eqref{A.18}. Then
\be
 F_l^{(k)}\ge F_{l,0}^{(k)}>\overline{F}_l^{(k)}\ge0. \label{A.30}
\ee
\end{lemma}

{\it Proof.}   Again we proceed by induction on $n\in \mathbb{N}$, such that  $1\le l-k\le n$. Hence for $n=1$, by direct inspection we easily check that \eqref{A.30} holds for $k=l-1$.
Hence we let $n>1,$ and by induction we assume that,
\be
 F_m^{(s)}\ge  F_{m,0}^{(s)}> \overline{F}_m^{(s)}\ge 0 \quad \forall\,m,s\in \mathbb{N}: \, 1\le m-s\le n-1,\label{A.31}
 \ee
and we are left to show that \eqref{A.30} holds also when  $l\ge 3$ and
\be
l-k=n.\label{A.32}
\ee

From \eqref{A.31} and \eqref{A.32} we see that,
\be
 F_l^{(j)}>0 \quad  F_{j-1}^{(k)}>0, \quad \forall\, k+1\le j\le l,\label{A.33}
\ee
 so we can use Proposition \ref{ap1} to conclude that $ F_l^{(k)}\ge  F_{l,0}^{(k)}$.

 To establish the second (strict) inequality in \eqref{A.30}, we use the derivation formulae \eqref{A.9}--\eqref{A.10}.   In view of \eqref{A.33} we find:
 \[
  \frac{\partial F_{l,0}^{(k)}}{\partial\alpha_{k,2}}=-\alpha_{k+1,1} F_l^{(k+2)} F_{k-1}^{(k)}=-\alpha_{k+1, 1} F_l^{(k+2)}<0
 \]
 so that,
 \[
 F_{l,0}^{(k)}(\alpha_{k,2}, \alpha_{k+1,1}, \dots, \alpha_{l,1})> F_{l,0}^{(k)}(1-\tau_k, \alpha_{k+1,1}, \dots, \alpha_{l,1}).
 \]
Furthermore,
\[\frac{\partial}{\partial\alpha_{k+1,1}}\left(F_{l,0}^{(k)}(1-\tau_k, \alpha_{k+1,1}, \dots, \alpha_{l,1})\right)=-(1-\tau_k)F_l^{(k+2)}<0,\]
which implies,
\[
 F_{l,0}^{(k)}(1-\tau_k, \alpha_{k+1,1}, \alpha_{k+1,2}\dots, \alpha_{l,1})> F_{l,0}^{(k)}(1-\tau_k, \tau_{k+1}, \alpha_{k+1,2} \dots, \alpha_{l,1}).
 \]

Thus, by observing that, for $k<j\le l$ we have:
\berr
&&\frac{\partial}{\partial\alpha_{j,2}}\left(F_{l,0}^{(k)}(1-\tau_k, \tau_{k+1}, 1-\tau_{k+1}, \dots, \tau_j,\alpha_{j,2}, \dots,  \alpha_{l,1})\right)\\
 &&=\frac{\partial F_{l,0}^{(k)}}{\partial\alpha_{j,2}}\left(1-\tau_k, \tau_{k+1}, 1-\tau_{k+1}, \dots, \tau_j,\alpha_{j,2}, \dots,  \alpha_{l,1}\right)\\
 &&=-\alpha_{j+1,1}F_l^{(j+2)}\overline{F}_{j-1}^{(k)}<0
\eerr
then we can proceed  inductively as above, to  conclude that,
\berr
 &&F_{l,0}^{(k)}(\alpha_{k,2}, \alpha_{k+1,1}, \alpha_{k+1,2}, \dots, \alpha_{l,1})>F_{l,0}^{(k)}(1-\tau_k, \alpha_{k+1,1}, \alpha_{k+1,2}, \dots, \alpha_{l,1})\\
 &&>F_{l,0}^{(k)}(1-\tau_k, \tau_{k+1}, \alpha_{k+1,2}, \dots, \alpha_{l,1})>F_{l,0}^{(k)}(1-\tau_k, \tau_{k+1}, 1-\tau_{k+1}, \dots, \alpha_{l,1})\\
&&>\cdots>F_l^{(k)}(1-\tau_k, \tau_{k+1}, 1-\tau_{k+1}, \dots, \tau_{l-1},1-\tau_{l-1}, \tau_l)=\overline{F}_l^{(k)}(\tau_k, \dots, \tau_l)\ge0
\eerr
and \eqref{A.30}  is established.

\hfill $\square$\\

{ \it Proof of Theorem \ref{ath1}}.  The property \eqref{A.19}  is a direct consequence of Lemma \ref{al4}.

\hfill $\square$

\begin{remark}
 By the tri-diagonal structure of the matrix $T_N^{(1)}$,  Lemma \ref{al4} and Sylvester's theorem, we can conclude that $T_N^{(1)}$ is actually
 positive definite, when \eqref{A.7}, \eqref{A.8} and \eqref{A.9} hold.
\end{remark}

\small{

}
\end{document}